\documentclass[12pt]{amsart}
\usepackage[]{amscd,amsfonts,amsmath,amsxtra,amssymb}
\usepackage{graphics}
\usepackage[all]{xy}
\usepackage{hyperref}
\usepackage[french]{babel}
\usepackage[T1]{fontenc}
\newtheorem{sub}{}[section]
\newtheorem{subsub}{}[sub]
\def\ov#1{\overline{#1}}

\def\coker{\mathop{\rm coker}\nolimits}
\def\Hom{\mathop{\rm Hom}\nolimits}
\def\HHom{\mathop{\mathcal Hom}\nolimits}
\def\OPL{\mathop\oplus\limits}
\def\Ext{\mathop{\rm Ext}\nolimits}
\def\EExt{\mathop{\mathcal Ext}\nolimits}
\def\Tor{\mathop{\rm Tor}\nolimits}

\def\Pic{\mathop{\rm Pic}\nolimits}

\def\End{\mathop{\rm End}\nolimits}
\def\EEnd{\mathop{\mathcal End}\nolimits}
\def\GL{\mathop{\rm GL}\nolimits}

\def\imm{\mathop{\rm im}\nolimits}
\def\deg{\mathop{\rm deg}\nolimits}
\def\rg{\mathop{\rm rg}\nolimits}
\def\Deg{\mathop{\rm Deg}\nolimits}
\def\spec{\mathop{\rm spec}\nolimits}
\def\lra{\longrightarrow}

\def\sigg{\mathop{\hbox{$\displaystyle\sum$}}\limits}

\def\hfl#1#2{\smash{\mathop{\ \hbox to 12mm{\rightarrowfill}}
\limits^{\scriptstyle#1}_{\scriptstyle#2} \ }}
\def\hflb#1#2{\smash{\mathop{\hbox to 12mm{\leftarrowfill}}
\limits^{\scriptstyle#1}_{\scriptstyle#2}}}

\def\m#1{{\hbox{$#1$}}}

\def\ot{\otimes}
\def\og{\leavevmode\raise.3ex\hbox{$\scriptscriptstyle\langle\!\langle$}}
\def\fg{\leavevmode\raise.3ex\hbox{$\scriptscriptstyle\,\rangle\!\rangle$}}

\def\nsp{\lbrace 0\rbrace}

\def\Ssect#1#2{\pagebreak[3]\begin{sub}\label{#2}{\sc\small\small 
#1}\rm\medskip}

\def\sepsec{\vskip 2.5cm}
\def\sepsub{\vskip 1.5cm}
\def\sepsubsub{\vskip 1cm}
\def\sepprop{\vskip 0.8cm}

\def\xmat#1{\[\xymatrix{#1}\]}
\def\flinc{\ar@{^{(}->}}
\def\fleq{\ar@{=}}
\def\flon{\ar@{->>}}
\def\fmaps{\ar@{|-{>}}}
\def\Nligne{\hfil\break}

\newcommand{\B}{{\mathbb B}}

\newcommand{\C}{{\mathbb C}}

\renewcommand{\P}{{\mathbb P}}
\newcommand{\D}{{\mathbb D}}
\newcommand{\F}{{\mathbb F}}
\newcommand{\E}{{\mathbb E}}
\newcommand{\G}{{\mathbb G}}

\newcommand{\U}{{\mathbb U}}
\newcommand{\V}{{\mathbb V}}
\newcommand{\W}{{\mathbb W}}
\renewcommand{\L}{{\mathbb L}}

\def\T{{\mathbb T}}

\newcommand{\kd}{{\mathcal D}}
\newcommand{\ke}{{\mathcal E}}
\newcommand{\kf}{{\mathcal F}}
\newcommand{\kg}{{\mathcal G}}

\newcommand{\ki}{{\mathcal I}}

\newcommand{\kk}{{\mathcal K}}
\newcommand{\kl}{{\mathcal L}}
\newcommand{\km}{{\mathcal M}}
\newcommand{\kn}{{\mathcal N}}
\newcommand{\ko}{{\mathcal O}}
\newcommand{\kp}{{\mathcal P}}

\newcommand{\ks}{{\mathcal S}}

\newcommand{\ku}{{\mathcal U}}
\newcommand{\kv}{{\mathcal V}}

\newcommand{\ky}{{\mathcal Y}}

\begin{document}

\def\refname{R\'ef\'erences}
\def\contentsname{Sommaire}
\def\proofname{D\'emonstration}
\def\abstractname{R\'esum\'e}

\author{Jean--Marc Dr\'{e}zet}
\address{
Institut de Math\'ematiques de Jussieu,
Case 247,
4 place Jussieu,
F-75252 Paris, France}
\email{drezet@math.jussieu.fr}
\urladdr{http://www.math.jussieu.fr/$\sim$drezet}

\title[{\tiny Faisceaux sans torsion et faisceaux quasi localement libres}]
{Faisceaux sans torsion et faisceaux quasi localement libres sur les courbes
multiples primitives}

\begin{abstract}
This paper is devoted to the study of some coherent sheaves on non reduced
curves that can be locally embedded in smooth surfaces. If $Y$ is such a curve
and $n$ is its multiplicity, then there is a filtration \m{C_1=C\subset
C_2\subset\cdots\subset C_n=Y} such that $C$ is the reduced curve associated to
$Y$, and for every \m{P\in C}, if \m{z\in\ko_{Y,P}} is an equation of $C$ then
\m{(z^i)} is the ideal of \m{C_i} in \m{\ko_{Y,P}}. A coherent sheaf on $Y$ is
called {\em torsion free} if it does not have any non zero subsheaf with finite
support. We prove that torsion free sheaves are reflexive. We study then the
{\em quasi locally free sheaves}, i.e. sheaves which are locally isomorphic to
direct sums of the \m{\ko_{C_i}}. We define an invariant for these sheaves, the
{\em complete type}, and prove the irreducibility of the set of sheaves of given
complete type. We study the generic quasi locally free sheaves, with
applications to the moduli spaces of stable sheaves on $Y$.
\end{abstract}

\maketitle

\tableofcontents

\section{Introduction}
\label{intro}

Une {\em courbe multiple primitive} est une vari\'et\'e alg\'ebrique complexe de
Cohen-Macaulay
pouvant localement \^etre plong\'ee dans une surface lisse, et dont la
sous-vari\'et\'e r\'eduite associ\'ee est une courbe lisse.
Les courbes projectives multiples primitives ont \'et\'e d\'efinies et
\'etudi\'ees pour la premi\`ere fois par C.~B\u anic\u a et O.~Forster dans
\cite{ba_fo}. Leur classification a \'et\'e faite dans \cite{ba_ei} pour les
courbes doubles, et dans \cite{dr1} dans le cas g\'en\'eral.

L'article \cite{dr2} donne les premi\`eres propri\'et\'es des faisceaux
coh\'erents et de leurs vari\'et\'es de modules sur les courbes multiples
primitives. Cette \'etude est poursuivie ici. En particulier la plupart des
r\'esultats de \cite{dr2} concernant les faisceaux sur les courbes doubles sont
g\'en\'eralis\'es ici en multiplicit\'e quelconque.

Les faisceaux semi-stables sur des vari\'et\'es non lisses ont d\'ej\`a \'et\'e
\'etudi\'es :
\begin{enumerate}
\item[--] sur des courbes r\'eduites dans \cite{ses}, \cite{bho}, \cite{bho2},
\cite{tei1}, \cite{tei2}, \cite{tei3}.
\item[--] sur des vari\'et\'es non r\'eduites semblables \`a celles qui sont
consid\'er\'ees ici dans\cite{in}.
\item[--] sur une vari\'et\'e ayant deux composantes irr\'eductibles
qui se coupent dans \cite{in2}.
\end{enumerate}
On peut esp\'erer en trouver des applications concernant les fibr\'es
vectoriels ou leurs vari\'et\'es de modules sur les courbes lisses
(cf. \cite{ei_gr}, \cite{sun0}, \cite{sun1}) en faisant d\'eg\'en\'erer des
courbes lisses vers une courbe multiple primitive. Le probl\`eme de la
d\'eg\'en\'eration des courbes lisses en courbes primitives doubles est
\'evoqu\'e dans \cite{gon}.

Les faisceaux coh\'erents sur les courbes non r\'eduites apparaissent aussi
lorsqu'on veut \'etudier les faisceaux de dimension 1 sur les surfaces,
lesquels interviennent notamment dans la conjecture de ``dualit\'e
\'etrange'' de J. Le Potier (cf. \cite{dan}).
Les faisceaux sur des courbes non r\'eduites apparaissent
(cf. \cite{lp}, \cite{lp4}) comme limites de fibr\'es vectoriels sur des
courbes lisses. Leur r\^ole est sans doute encore plus important si on cherche
\`a obtenir d'autres vari\'et\'es de modules fins de faisceaux de dimension 1
que les classiques vari\'et\'es de modules de faisceaux semi-stables (cf.
\cite{dr}).

\sepsubsub

{\bf Notations :} Si $\ke$, $\kf$ sont des faisceaux coh\'erents sur une
vari\'et\'e alg\'ebrique, on d\'esigne par \m{\Ext^i(\ke,\kf)} (resp.
\m{\Hom(\ke,\kf)}, \m{\End(\ke)}) les Ext (resp. morphismes) globaux de
faisceaux et par \m{\EExt^i(\ke,\kf)} (resp. \m{\HHom(\ke,\kf)}, \m{\EEnd(\ke)})
les faisceaux d'Ext locaux (resp. de morphismes).

\sepsub

\Ssect{Faisceaux sans torsion}{fst}

Soit $Y$ une courbe multiple primitive.
Un faisceau coh\'erent sur $Y$ est dit {\em sans torsion} s'il n'a pas de
sous-faisceau non nul de support fini. Les faisceaux semi-stables sur $Y$
non concentr\'es en un nombre fini de points sont des exemples de faisceaux sans
torsion. On d\'emontre en \ref{theo_refl} le

\sepprop

{\bf Th\'eor\`eme A : }{\em Soit $\ke$ un faisceau
coh\'erent sur $Y$. Alors les propri\'et\'es suivantes sont \'equivalentes:
\begin{itemize}
\item[(i)] $\ke$ est r\'eflexif.
\item[(ii)] $\ke$ est sans torsion.
\item[(iii)] On a \ $\EExt^1_{\ko_Y}(\ke,\ko_Y)=0$ .
\end{itemize}
Si les conditions pr\'ec\'edentes sont r\'ealis\'ees on a de plus \
\m{\EExt^i_{\ko_Y}(\ke,\ko_Y)=0} \ pour tout \m{i\geq 1}.
}

\sepprop

Ce r\'esultat avait \'et\'e d\'emontr\'e pour les courbes doubles dans
\cite{dr2}.
Ceci implique en particulier qu'en dualisant une suite exacte de faisceaux sans
torsion sur \m{C_n} on obtient une suite exacte. Une autre application est le
{\em th\'eor\`eme de dualit\'e de Serre} sur $Y$, qui est vrai pour les
faisceaux sans torsion : si $\ke$ un faisceau
coh\'erent sans torsion sur $Y$, alors on a des isomorphismes fonctoriels
\[H^i(Y,\kf) \ \simeq \ H^{1-i}(Y,\kf^\vee\ot\omega_Y)^*\]
pour \m{i=0,1}.

\end{sub}

\sepsub

\Ssect{Faisceaux quasi localement libres}{fql}

Soit $Y$ une courbe multiple primitive de courbe r\'eduite associ\'ee $C$
projective et de multiplicit\'e $n$. Soient \m{\ki_C} le faisceau d'id\'eaux de
$C$, et pour \m{1\leq i\leq n}, \m{C_i} le sous-sch\'ema de $Y$ d\'efini par
\m{\ki_C^i}. On a donc \m{C_n=Y} et \m{C_i} est une courbe multiple primitive de
multiplicit\'e $i$. On note \ \m{L=\ki_C/\ki_C^2} , qui est un fibr\'e en
droites sur $C$.

Si $\kf$ est un faisceau coh\'erent sur $Y$ on note \m{\kf_i} le noyau de la
restriction \ \m{\kf\to\kf_{\mid C_i}}. Les quotients
\ \m{G_i(\kf) \ = \ \kf_i/\kf_{i+1}} ,
\m{0\leq i<n}, sont des faisceaux sur $C$. Ils permettent de d\'efinir les {\em
rang g\'en\'eralis\'e} et le {\em degr\'e g\'en\'eralis\'e} de $\kf$ :
\[R(\kf)=\sigg_{i=0}^{n-1}\rg(G_i(\kf)) , \quad\quad
\Deg(\kf)=\sigg_{i=0}^{n-1}\deg(G_i(\kf))\]
(cf. \ref{inva}, \cite{dr2})

Un faisceau coh\'erent $\ke$ sur $Y$ est dit {\em quasi localement libre} s'il
est localement isomorphe \`a une somme directe du type
\[\bigoplus_{i=1}^{n}m_i\ko_{C_i} \ .\]
La suite \m{(m_1,\ldots,m_{n})} s'appelle le {\em type} de $\ke$. D'apr\`es
\cite{dr2} pour que $\ke$ soit quasi localement libre il faut et il suffit que
tous les \m{G_i(\ke)} soient localement libres sur $C$. La paire
\[\big((\rg(G_0(\ke)),\ldots,\rg(G_{n-1}(\ke)) \ , \
(\deg(G_0(\ke)),\ldots,\deg(G_{n-1}(\ke))\big)\]
s'appelle le {\em type complet} de $\ke$.

\sepprop

\begin{subsub}Invariance du type complet par d\'eformation - \rm Soient
\m{m_1,\ldots,m_{n}\geq 0} des entiers, $Y$ une vari\'et\'e alg\'ebrique
int\`egre et \m{\kf} un faisceau coh\'erent sur \m{Y\times C_n}, plat sur $Y$.
On suppose que pour tout point ferm\'e \m{y\in Y}, \m{\kf_y} est quasi
localement libre de type \m{(m_1,\ldots,m_{n})}. Soit \ \m{p_{C_n}:Y\times
C_n\to C_n} \ la projection. On d\'emontre en \ref{f_ql_th1} le

\sepprop

{\bf Th\'eor\`eme B : }{\em Pour tout point ferm\'e \m{(y,P)\in Y\times C_n} il
existe un voisinage $U$ de \m{(y,P)} tel que
\[\kf_{\mid U} \ \simeq \ \bigoplus_{i=1}^{n}m_i.p_{C_n}^*(\ko_{C_i})_{\mid
 U}.\]}

\sepprop

Il en d\'ecoule que le noyau \m{\kf_i} de la restriction \ \m{\kf\to\kf_{\mid
Y\times C_i}} \ est plat sur $Y$ et que pour tout \m{y\in Y} on a \
\m{(\kf_i)_y=(\kf_y)_i}. On en d\'eduit ais\'ement que le type complet de
\m{\kf_y} est ind\'ependant de \m{y\in Y}.
\end{subsub}

\sepprop

\begin{subsub} Irr\'eductibilit\'e - \rm Soient $X$ une vari\'et\'e
alg\'ebrique, et $\ky$ un ensemble de classes d'isomorphisme de faisceaux
coh\'erents sur $X$. On dit que $\ky$ est {\em irr\'eductible} si pour tous
\m{E_0,E_1\in\ky} il existe une vari\'et\'e irr\'eductible $Z$, un faisceau
coh\'erent $\ke$ sur \m{Z\times X} plat sur $Z$, tel que pour tout point ferm\'e
\m{z\in Z} on ait \m{\ke_z\in\ky}, et qu'il existe des points ferm\'es
\m{z_0,z_1\in Z} tels que \m{\ke_{z_0}=E_0}, \m{\ke_{z_1}=E_1}. Il est bien
connu par exemple que les fibr\'es vectoriels de rang et degr\'e donn\'es sur
$C$ constituent un ensemble irr\'eductible. On d\'emontre en \ref{irred_theo} le

\sepprop

{\bf Th\'eor\`eme C : }{\em Les faisceaux quasi localement libres sur \m{C_n} de
type complet donn\'e constituent un ensemble irr\'eductible.}
\end{subsub}

\sepprop

\begin{subsub}\label{intro3}Faisceaux quasi localement libres de type rigide -
\em Un faisceau quasi localement libre est dit {\em de type rigide} s'il est
localement libre, ou localement isomorphe \`a \m{a\ko_{C_n}\oplus\ko_{C_k}},
avec \m{1\leq k<n}. On d\'emontre en \ref{qlltr1} la

\sepprop

{\bf Proposition D : }{\em Soient $Y$ une vari\'et\'e
alg\'ebrique int\`egre et $\kf$ une famille plate de faisceaux coh\'erents sur
\m{C_n} param\'etr\'ee par $Y$. Alors l'ensemble des points \m{y\in Y} tels que
\m{\ke_y} soit quasi localement libre de type rigide est un ouvert de $Y$.}

\sepprop

Les faisceaux quasi localement libres de type rigide sont donc des faisceaux
g\'en\'eriques.

Soit $\ke$ un faisceau coh\'erent sur \m{C_n}. Soit
\m{(\widetilde{\ke},S,s_0,\epsilon)} une {\em d\'eformation semi-universelle}
de $\ke$ (cf. \cite{si_tr}, \cite{dr3} 3.1), donc $\widetilde{\ke}$ est une
famille plate de faisceaux coh\'erents sur \m{C_n} param\'etr\'ee par $S$,
\m{s_0} est un point ferm\'e de $S$ et \
\m{\epsilon:\widetilde{\ke}_{s_0}\simeq\ke}. Le morphisme de d\'eformation
infinit\'esimale de Koda\"ira-Spencer
\[\omega_{\widetilde{\ke},s_0}:T_{s_0}S\lra\Ext^1_{\ko_{C_n}}(\ke,\ke)\]
est un isomorphisme. On pose
\[D_{reg}(\ke) \ = \ \omega_{\widetilde{\ke},s_0}(T_{s_0}(S_{red})) .\]
Pour toute d\'eformation de $\ke$ param\'etr\'ee par une vari\'et\'e r\'eduite,
le morphisme de d\'eformation infinit\'esimale de Koda\"ira-Spencer
correspondant est \`a valeurs dans \m{D_{reg}(\ke)}.

On dit que $\ke$ est {\em lisse} si $S$ est lisse en \m{s_0}. Les faisceaux
localement libres sont lisses mais ce n'est pas le cas en g\'en\'eral des
autres faisceaux quasi localement libres de type rigide.

On dit que $\ke$ est {\em lisse pour les d\'eformations r\'eduites} si
\m{S_{red}} est lisse en \m{s_0}.

On a une suite exacte canonique
\[0\lra H^1(\EEnd(\ke))\lra\Ext^1_{\ko_{C_n}}(\ke,\ke)\lra
H^0(\EExt^1_{\ko_{C_n}}(\ke,\ke))\lra 0 .\]
On d\'emontre en \ref{qlltr2} le

\sepprop

{\bf Th\'eor\`eme E : }{\em Si $\ke$ est un faisceau quasi
localement libre de type rigide g\'en\'erique, alors on a \ \m{D_{reg}(\ke)=
H^1(\EEnd(\ke))}.}

\sepprop

d'o\`u on d\'eduit le

\sepprop

{\bf Corollaire F : }{\em Si $\ke$ est un faisceau
quasi localement libre de type rigide tel que \m{\dim_\C(\End(\ke))} soit
minimal, c'est \`a-dire soit tel que pour tout faisceau quasi localement libre
$\kf$ de m\^eme type complet que $\ke$ on ait \ \m{\dim_\C(\End(\kf))\geq
\dim_\C(\End(\ke))}, alors on a \ \m{D_{reg}(\ke)= H^1(\EEnd(\ke))} , et $\ke$
est lisse pour les d\'eformations r\'eduites.}
\end{subsub}

\sepprop

\begin{subsub} Composantes des vari\'et\'es de modules de faisceaux sans
torsion stables -
\rm Soient $R$, $D$ des entiers, avec \m{R\geq 1}. Soit \m{\km(R,D)} la
vari\'et\'e de modules des faisceaux stables de rang g\'en\'eralis\'e $R$
et de degr\'e g\'en\'eralis\'e $D$ sur \m{C_n} (cf. \cite{ma1}, \cite{ma2},
\cite{si}). On supposera que \m{\deg(L)<0}, car dans le cas contraire les
seuls faisceaux sans torsion stables sur \m{C_n} sont les fibr\'es vectoriels
stables sur $C$.

Soit $\ke$ un faisceau quasi localement
libre de type rigide localement isomorphe \`a \ \m{a\ko_n\oplus\ko_k} . Posons
\[E=\ke_{\mid C} , \quad F=G_k(\ke)\ot L^{-k} , \quad \delta(\ke)=\deg(F) ,
\quad \epsilon(\ke)=\deg(E) .\]
On a
\[R(\ke)=an+k , \quad \Deg(\ke) = k\epsilon(\ke)+(n-k)\delta(\ke)
+\big(n(n-1)a+k(k-1)\big)\deg(L)/2 \ .\]
D'apr\`es \ref{intro3} les d\'eformations de $\ke$ sont des
faisceaux quasi localement libres de type rigide, et \m{a(\ke)}, \m{k(\ke)},
\m{\delta(\ke)} et \m{\epsilon(\ke)} sont aussi invariants par d\'eformation.

Soient \ \m{R=R(\ke), D=\Deg(\ke), \delta=\delta(\ke), \epsilon=\epsilon(\ke)}.
Les faisceaux quasi localement libres de type rigide stables $\kf$ tels que
\m{a(\kf)=a}, \m{k(\kf)=k}, \m{\delta(\kf)=\delta}, \m{\epsilon(\kf)=\epsilon}
constituent donc un ouvert de \m{\km(R,D)}, not\'e \m{\kn(a,k,\delta,\epsilon)}.

\sepprop

{\bf Proposition G : }{\em La vari\'et\'e \m{\kn(a,k,\delta,\epsilon)} est
irr\'eductible, et la sous-vari\'et\'e r\'eduite sous-jacente
\m{\kn(a,k,\delta,\epsilon)_{red}} est lisse. Si cette vari\'et\'e est non
vide, on a
\[\dim(\kn(a,k,\delta,\epsilon)) \ = \ 1 - \big(\frac{n(n-1)}{2}a^2+k(n-1)a+
\frac{k(k-1)}{2}\big)\deg(L)+(g-1)(na^2+k(2a+1))\]
($g$ d\'esignant le genre de $C$).
Pour tout faisceau $\kf$ de \m{\kn(a,k,\delta,\epsilon)_{red}}, l'espace
tangent de \m{\kn(a,k,\delta,\epsilon)_{red}} en $\kf$ est canoniquement
isomorphe \`a \m{H^1(\EEnd(\kf))} .}

\sepprop

Les param\`etres \m{a,k,\delta,\epsilon} tels que \m{\kn(a,k,\delta,\epsilon)}
soit non vide ne sont pas encore connus. Des exemples pour \m{a=1}, \m{k=1} sur
une courbe double sont donn\'es dans \cite{dr2}.
\end{subsub}

\end{sub}

\sepsub

\Ssect{Questions}{quest}

\begin{subsub}Composantes irr\'eductibles - \rm Soit \m{\kk(R,D)} l'ensemble
des classes d'isomorphisme de faisceaux sans torsion de rang g\'en\'eralis\'e
$R$ et de degr\'e g\'en\'eralis\'e $D$ sur \m{C_n}. Il serait int\'eressant de
d\'ecomposer \m{\kk(R,D)} en sous-ensembles irr\'eductibles et de param\'etrer
ces ``composantes''. On a d\'ecrit dans le pr\'esent article celles qui
contiennent des faisceaux quasi localement libres de type rigide. Mais pour
certaines valeurs de $R$ et $D$ de telles composantes n'existent pas. Supposons
par exemple que $R$ soit multiple de $n$ : \m{R=an}. Si $\ke$ est quasi
localement libre de type rigide et de rang g\'en\'eralis\'e $R$, alors $\ke$
est localement libre et on a
\[\Deg(\ke) \ = \ n.\deg(\ke_{\mid C})+\frac{n(n-1)}{2}a\deg(L) ,\]
Donc si \ \m{D\not\equiv\frac{n(n-1)}{2}a\deg(L)\ ({\rm mod}\ n)} \ il n'existe
pas de faisceau quasi localement libre de type rigide de rang $R$ et de degr\'e
g\'en\'eralis\'e $D$. Les faisceaux g\'en\'eriques semblent \^etre dans de cas
des faisceaux localement isomorphes \`a \m{a\ko_{C_n}} sur un ouvert non vide de
\m{C_n}, mais avec un certain nombre de points singuliers.
\end{subsub}

\sepprop

\begin{subsub}Conditions d'existence des faisceaux stables - \rm Quels sont les
entiers $R$, $D$ tels qu'il existe un faisceau stable de rang $R$ et de
degr\'e $D$ sur \m{C_n} ? On doit supposer que \m{\deg(L)<0}, sans quoi il
n'existe aucun faisceau stable en dehors des fibr\'es vectoriels stables sur
$C$. Autre question analogue : Quels sont les entiers $R$, $D$ tels qu'il existe
un faisceau simple de rang $R$ et de degr\'e $D$ sur \m{C_n} ?
\end{subsub}

\end{sub}

\sepsub

\Ssect{Plan des chapitres suivants}{PlanC}

Le chapitre 2 contient des r\'esultats techniques utilis\'es dans les autres
chapitres.

Dans le chapitre 3 on rappelle en les approfondissant certaines propri\'et\'es
des faisceaux coh\'erents sur les courbes multiples primitives. En particulier
on \'etudie des m\'ethodes de construction de faisceaux coh\'erents sur \m{C_n}
qui sont utilis\'ees dans les chapitres suivants. Par exemple on donne la
construction des faisceaux $\ke$ connaissant \m{\ke_{\mid C}} et le noyau de
\m{\ke\to\ke_{\mid C}}. Ces m\'ethodes sont utilis\'ees dans des
d\'emonstrations par r\'ecurrence sur la multiplicit\'e $n$.

Le chapitre 4 traite des faisceaux r\'eflexifs. Les r\'esultats les plus
importants d\'emontr\'es ici sont le th\'eor\`eme A et la dualit\'e de Serre
pour les faisceaux sans torsion.

Le chapitre 5 contient un r\'esulat technique indispensable aux d\'emonstrations
du chapitre 6~: si $\ke$ est un faisceau quasi localement libre, alors il existe
un fibr\'e vectoriel $\E$ sur \m{C_n} et un morphisme surjectif \m{\E\to\ke}
induisant un isomorphisme \m{\E_{\mid C}\simeq\ke_{\mid C}}. On en d\'eduit
l'existence de curieuses r\'esolutions localement libres des faisceaux quasi
localement libres (les {\em r\'esolutions p\'eriodiques}), qui permettent de
mettre en \'evidence une classe canonique
\[\lambda_\ke \ \in \ \Ext^2_{\ko_{C_n}}(\ke,\ke)\]
associ\'ee \`a tout faisceau quasi localement libre $\ke$ sur \m{C_n}, qui est
nulle si et seulement si $\ke$ est localement libre.

Le chapitre 6 contient les d\'emonstrations de r\'esultats \'enonc\'es en
\ref{fql}.
\end{sub}

\sepsec
\newpage
\section{Pr\'eliminaires}\label{prelim}

\Ssect{D\'efinitions des courbes multiples primitives et notations}{def_nota}

Une {\em courbe primitive} est une vari\'et\'e $Y$ de
Cohen-Macaulay telle que la sous-vari\'et\'e r\'eduite associ\'ee \m{C=Y_{red}}
soit une courbe lisse irr\'eductible, et que tout point ferm\'e de
$Y$ poss\`ede un voisinage pouvant \^etre plong\'e dans une surface lisse.

Soient $P$ un point ferm\'e de $Y$, et $U$ un voisinage de $P$ pouvant
\^etre plong\'e dans une surface lisse $S$. Soit $z$ un \'el\'ement de
l'id\'eal maximal de l'anneau local \m{\ko_{S,P}} de $S$ en $P$ engendrant
l'id\'eal de $C$ dans cet anneau. Il existe alors un unique entier $n$,
ind\'ependant de $P$, tel que l'id\'eal de $Y$ dans \m{\ko_{S,P}} soit
engendr\'e par \m{(z^n)}. Cet entier $n$ s'appelle la {\em multiplicit\'e} de
$Y$. Si \m{n=2} on dit que $Y$ est une {\em courbe double}. D'apr\`es
\cite{dr1}, th\'eor\`eme 5.2.1, l'anneau \m{\ko_{Y,P}} est isomorphe \`a
\m{\ko_{CP}\ot(\C[t]/(t^n))}.

Soit \m{\ki_C} le faisceau d'id\'eaux de $C$ dans $Y$. Alors le faisceau
conormal de $C$, \m{L=\ki_C/\ki_C^2} est un fibr\'e en droites sur $C$, dit {\em
associ\'e} \`a $Y$. Il existe une filtration canonique
\[C=C_1\subset\cdots\subset C_n=Y\ ,\]
o\`u au voisinage de chaque point $P$ l'id\'eal de \m{C_i} dans \m{\ko_{S,P}}
est \m{(z^i)}. On notera \ \m{\ko_i=\ko_{C_i}} \ pour \m{1\leq i\leq n}.

Le faisceau \m{\ki_C} est un fibr\'e en droites sur \m{C_{n-1}}. Il existe
d'apr\`es \cite{dr2}, th\'eor\`eme 3.1.1, un fibr\'e en droites $\L$ sur
\m{C_n} dont la restriction \`a \m{C_{n-1}} est \m{\ki_C}. On a alors, pour
tout faisceau de \m{\ko_n}-modules $\ke$ un morphisme canonique
\[\ke\ot\L\lra\ke\]
qui en chaque point ferm\'e $P$ de $C$ est la multiplication par $z$.

Pour tout point ferm\'e $P$ de $C$, on notera \m{z_P\in\ko_{nP}} une \'equation
de $C$ et \m{x_P\in\ko_{nP}} un \'el\'ement au dessus d'un g\'en\'erateur de
l'id\'eal maximal de \m{\ko_{CP}}. S'il n'y a pas d'ambigu\"it\'e on notera $z$,
$x$ au lieu de \m{z_P}, \m{x_P} respectivement.

Soient $X$, $Y$, $T$ des vari\'et\'es alg\'ebriques. On note \m{p_X}, \m{p_Y}
les projections \m{X\times Y\to X}, \m{X\times Y\to Y}.

Si $\ke$ est un faisceau coh\'erent sur \m{Y\times T} et si \m{f:X\to Y} est un
morphisme, on notera \ \m{f^\sharp(\ke)=(f\times I_T)^*(\ke)} .
\end{sub}

\sepsub

\Ssect{Fibr\'e canonique}{fibcan}

Soit \m{C_n} une courbe multiple primitive de multiplicit\'e $n$ de courbe
r\'eduite
associ\'ee $C$. On suppose que $C$ est projective, il en est donc de m\^eme de
\m{C_n}. \'Etant localement intersection compl\`ete, \m{C_n} admet un {\em
faisceau dualisant} \m{\omega_{C_n}}, qui est un fibr\'e en droites sur \m{C_n}.
On peut le d\'efinir \`a partir d'un plongement de \m{C_n} dans une vari\'et\'e
projective lisse $X$ :
\[\omega_{C_n} \ \simeq \ \omega_X\ot\det(\ki/\ki^2)^* ,\]
$\ki$ d\'esignant le faisceau d'id\'eaux de \m{C_n} dans $X$ (on peut voir
\m{\ki/\ki^2} comme un fibr\'e vectoriel de rang \m{\dim(X)-1} sur \m{C_n}).

\end{sub}

\sepsub

\Ssect{Les $\Ext$ de faisceaux d\'efinis sur des sous-vari\'et\'es}{extf}

Soient $X$ une vari\'et\'e projective et \m{Y\subset X} une sous-vari\'et\'e
ferm\'ee. Si \ \m{j:Y\to X} \ est l'inclusion et $E$ un faisceau coh\'erent sur
$Y$, on notera aussi souvent $E$ le faisceau \m{j_*(E)} sur $X$.

\sepprop

\begin{subsub}{\bf Proposition : }\label{extf_pr} Soient $E$, $F$ des faisceaux
coh\'erents sur $Y$.

1 - On a une suite exacte fonctorielle canonique
\[0\lra\Ext^1_{\ko_Y}(F,E)\lra\Ext^1_{\ko_X}(F,E)\lra\Hom(
\Tor^1_{\ko_X}(F,\ko_Y),E)\lra\Ext^2_{\ko_Y}(F,E)  .\]

2 - On a une suite exacte fonctorielle canonique
\[0\lra\EExt^1_{\ko_Y}(F,E)\lra\EExt^1_{\ko_X}(F,E)\lra\HHom(
\Tor^1_{\ko_X}(F,\ko_Y),E)\lra\EExt^2_{\ko_Y}(F,E)  .\]
\end{subsub}

\begin{proof}
La premi\`ere assertion est la proposition 2.2.1 de \cite{dr2}, et la seconde se
d\'emontre de mani\`ere analogue.
\end{proof}

\end{sub}

\sepsub

\Ssect{Extensions}{extens}

Soient $X$ une vari\'et\'e alg\'ebrique et $E$, $F$ des faisceaux coh\'erents
sur $X$. Une {\em extension} de $F$ par $E$ est une suite exacte de faisceaux
de \m{\ko_X}-modules \ \m{0\to E\to\ke\to F\to 0} . Il est bien connu que les
classes d'isomorphisme d'extensions sont param\'etr\'ees par
\m{\Ext^1_{\ko_X}(F,E)}. On a une suite exacte canonique
\xmat{0\ar[r] & H^1(\HHom(F,E))\ar[r]^-\pi & \Ext^1_{\ko_X}(F,E)\ar[r]^-p &
H^0(\EExt^1_{\ko_X}(F,E))\ar[r] & 0 ,}
d'apr\`es la suite spectrale des Ext (cf. \cite{go}, 7.3). Le morphisme $p$ est
facile \`a d\'ecrire : soient \m{x\in X}, \m{\sigma\in\Ext^1_{\ko_X}(F,E)} et \
\m{0\to E\to\ke\to F\to 0} \ l'extension associ\'ee. Alors \m{p(\sigma)(x)} est
l'\'el\'ement de \m{\EExt^1_{\ko_{Xx}}(F_x,E_x)} associ\'e \`a la suite exacte
\ \m{0\to E_x\to\ke_x\to F_x\to 0} . Les extensions associ\'ees aux
\'el\'ements de \m{\imm(\pi)} sont donc celles qui sont triviales en chaque
point de $X$. Pour les d\'ecrire plus pr\'ecis\'ement on va rappeler une
construction explicite des extensions.

\sepprop

\begin{subsub}\label{ext_cons} Construction des extensions - \rm
Soit \m{\ko_X(1)} un fibr\'e en droites tr\`es ample sur $X$. Si
\m{m_0\gg 0}, il existe un entier \m{k_0>0} et un morphisme surjectif
\[f_0:\ko_X(-m_0)\ot\C^{k_0}\lra F .\]
Si \m{m_1\gg 0} il existe de m\`eme un entier \m{k_1>0} et un morphisme
surjectif
\[f_1:\ko_X(-m_1)\ot\C^{k_1}\lra\ker(f_0) .\]
On obtient en poursuivant ainsi une r\'esolution localement libre de $F$ :
\xmat{\cdots\ar[r] & F_2\ar[r]^-{f_2} & F_1\ar[r]^-{f_1} & F_0\flon[r]^-{f_0} &
F \ ,}
avec \ \m{F_i=\ko_X(-m_i)\ot\C^{k_i}} . Soit \ \m{\phi_i:\Hom(F_{i-1},E)\to
\Hom(F_i,E)} \ le morphisme induit par \m{f_i}. Si \m{m_0,m_1,m_2\gg 0} on a un
isomorphisme canonique
\[\Ext^1_{\ko_X}(F,E) \ \simeq \ \ker(\phi_2)/\imm(\phi_1) .\]
Soient \ \m{\sigma\in\Ext^1_{\ko_X}(F,E)} \ et \m{\phi:F_1\to E} au dessus de
$\sigma$ (on a donc \m{\phi\circ f_2=0}). Soient \m{\eta=(\phi,f_1)} et
\m{\ke=\coker(\eta)}. On a un morphisme injectif \'evident \m{E\to\ke}, et un
morphisme surjectif \m{\ke\to F}, dont le noyau est \'egal \`a l'image du
pr\'ec\'edent, d'o\`u une suite exacte \ \m{0\to E\to\ke\to F\to 0} , dont
l'\'el\'ement associ\'e de \m{\Ext^1_{\ko_X}(F,E)} est pr\'ecis\'ement $\sigma$.

Un \'el\'ement $\lambda$ de l'image de $\pi$ provient d'un morphisme
\m{\psi:F_1\to E} se factorisant localement par \m{F_0}, c'est-\`a-dire qu'il
existe un recouvrement ouvert \m{(U_i)} de $X$ et pour tout $i$ un morphisme \
\m{\alpha_i:F_{0\mid U_i}\to E_{\mid U_i}} \ tel que sur \m{U_i} ont ait \
\m{\psi=\alpha_if_1}. Pour tous \m{i,j}, \m{\alpha_i-\alpha_j} s'annule sur
\m{\ker(f_0)=\imm(f_1)} et induit donc un morphisme \ \m{\tau_{ij}:F_{\mid
U_{ij}}\to E_{\mid U_{ij}}} . La famille \m{(\tau_{ij})} est un cocycle et
d\'efinit donc un \'el\'ement de \m{H^1(\HHom(F,E))}, dont l'image dans
\m{\Ext^1_{\ko_X}(F,E)} est $\lambda$.

R\'eciproquement, \'etant donn\'e un tel cocycle \m{(\tau_{ij})},
\m{(\tau_{ij}f_0)} repr\'esente un \'el\'ement de\hfil\break
\m{H^1(\HHom(F_0,E))}. On peut supposer ce dernier nul (car \m{m_0\gg 0}). Donc
il existe une famille \m{(\alpha_i)}, \m{\alpha_i:F_{0\mid U_i}\to E_{\mid U_i}}
, telle que pour tout $i$, $j$ on ait \ \m{\tau_{ij}f_0=\alpha_i-\alpha_j} . On
a \ \m{\alpha_if_1=\alpha_jf_1} \ sur \m{U_{ij}}, donc les \m{\alpha_if_1} se
recollent et d\'efinissent un morphisme \m{F_1\to E} d'o\`u provient
l'\'el\'ement de \m{\Ext^1_{\ko_X}(F,E)} induit par \m{(\tau_{ij})}.
\end{subsub}

\sepprop

\begin{subsub}\label{ext_p}Description de $\pi$ - \rm Soient \m{\sigma\in
\Ext^1_{\ko_X}(F,E)} et \m{\lambda\in H^1(\HHom(F,E))}. Soit \m{(U_i)} un
recouvrement ouvert de $X$ tel que $\lambda$ soit repr\'esent\'e
par un cocycle \m{(\tau_{ij})}, \m{\tau_{ij}:F_{\mid U_{ij}}\to E_{\mid
U_{ij}}}. Soit
\xmat{0\ar[r] & E\ar[r]^\nu & \ke\ar[r]^\mu & F\ar[r] & 0}
l'extension correspondant \`a $\sigma$. Posons
\[\gamma_{ij} = I_\ke+\nu\tau_{ij}\mu : \ke_{\mid U_{ij}}\lra\ke_{\mid U_{ij}}
.\]
Alors on a \ \m{\gamma_{ij}\gamma_{jk}=\gamma_{ik}} \ pour tous \m{i,j,k}. Donc
on peut d\'efinir un faisceau coh\'erent $\kf$ sur $X$ en recollant les
\m{\ke_{U_{ij}}} au moyen des automorphismes \m{\gamma_{ij}}. Les morphismes
\ \m{\nu:E_{\mid U_{ij}}\to\ke_{\mid U_{ij}}} \ et \ \m{\mu:\ke_{\mid U_{ij}}\to
F_{\mid U_{ij}}} \ se recollent aussi (car \m{\mu\gamma_{ij}=\mu} et
\m{\gamma_{ij}\nu=\nu}), et on obtient une suite exacte
\xmat{0\ar[r] & E\ar[r]^-{\nu'} & \kf\ar[r]^-{\mu'} & F\ar[r] & 0 ,}
correspondant \`a \m{\sigma'\in\Ext^1_{\ko_X}(F,E)}.
\end{subsub}

\sepprop

\begin{subsub}\label{ext_P} {\bf Proposition : } On a \
\m{\sigma'=\sigma+\lambda} . 
\end{subsub}

\begin{proof} On utilise les notations pr\'ec\'edentes. On a \
\m{\kf=\coker(\eta')}, o\`u
\[\eta'_{\mid U_i}=(\phi+\alpha_if_1,f_1):F_{1\mid U_i}\lra E_{\mid U_i}\oplus
 F_{0\mid U_i} .\]
Rappelons que \m{\ke=\coker(\eta)}, avec \m{\eta=(\phi,f_1)}. On a donc un
diagramme commutatif sur \m{U_i}
\xmat{F_1\ar[r]^-{\eta'}\fleq[d] & E\oplus F_0\ar[d]^{\beta_i} \\
F_1\ar[r]^-\eta & E\oplus F_0}
avec \ \m{\beta_i=\begin{pmatrix}1 & -g_i\\ 0 & 1\end{pmatrix}} , induisant un
isomorphisme \ \m{\gamma_i:\kf_{\mid U_i}\simeq\ke_{\mid U_i}} . Le r\'esultat
d\'ecoule du fait que \ \m{\gamma_j\gamma_i^{-1}=\gamma_{ij}} (cf. \ref{ext_p}).
\end{proof}

\end{sub}

\sepsub

\Ssect{Prolongement de fibr\'es vectoriels}{prol_f}

Soit \m{C_n} une courbe multiple primitive de multiplicit\'e $n$ de courbe
r\'eduite associ\'ee $C$. Le r\'esultat suivant g\'en\'eralise le th\'eor\`eme
3.1.1 de \cite{dr2} :

\sepprop

\begin{subsub}\label{prol_th}{\bf Th\'eor\`eme :} Soient $X$ une vari\'et\'e
alg\'ebrique affine irr\'eductible, $i$ un entier tel que \m{1\leq i<n} et $Y$
une sous-vari\'et\'e ferm\'ee de \m{X\times C_n} contenant \m{X\times C_i}.
Soit $\ke$ un faisceau localement libre sur $Y$. Alors il existe un faisceau
localement libre $\E$ sur \m{X\times C_n} tel que \ \m{\E_{\mid Y}\simeq\ke}.
\end{subsub}

\begin{proof}
On utilisera le r\'esultat suivant : pour tout faisceau coh\'erent $\kf$ sur
\m{X\times C_n}, on a \m{H^i(\kf)=\nsp} pour \m{i\geq 2}. En effet, soit
\m{p_n} la projection \m{X\times C_n\to C_n}. Comme $X$ est affine, les images
directes \m{R^jp_n^*(\ke)} sont nulles pour \m{j>0}. On a donc \ \m{H^i(\kf)
\simeq H^i(p_{n*}(\kf))}, qui est nul si \m{i\geq 2}.

En raisonnant par r\'ecurrence on se ram\`ene au cas o\`u \m{i=n-1}.

Soit \m{\ki_Y} le faisceau d'id\'eaux de $Y$ dans \m{X\times C_n}. Remarquons
que son support contenu dans \m{X\times C}. On note \m{E_0=\ke_{\mid X\times
C}}.

Il existe un recouvrement ouvert \m{(U_i)} de \m{Y} tel que chaque
restriction \m{\ke_{\mid U_i}} soit un fibr\'e trivial. On peut voir \m{(U_i)}
comme un recouvrement de \m{X\times C_n}. Soient \ \m{\lambda_i:\ke_{\mid
U_i}\simeq\ko_Y(U_i)\ot\C^r} \ des trivialisations, \
\m{\lambda_{ij}=\lambda_j\circ\lambda_i^{-1}\in\GL(r,\ko_Y(U_{ij}))}.
Soient \m{\Lambda_{ij}\in\GL(r,\ko_{X\times C_n}(U_{ij}))} une extension de
\m{\lambda_{ij}} et \ \m{\rho_{ijk}=\Lambda_{jk}\Lambda_{ij}-\Lambda_{ik}} \
(\'el\'ement de \m{\ko_{X\times C_n}(U_{ijk})\ot\End(\C^r)}). Alors les
\m{\Lambda_{ij}} d\'efinissent un fibr\'e vectoriel sur \m{X\times C_n} si et
seulement si les \m{\rho_{ijk}} sont nuls. Leurs restrictions \`a \m{Y}
sont nulles, donc on peut les consid\'erer comme des \'el\'ements de \
\m{\End(\ko_{X\times C_n}(U_{ijk})\ot\C^r)\ot\ki_Y(U_{ijk})}.
Soit \
\m{\mu_{ijk}=(\lambda_k)^{-1}\rho_{ijk}\lambda_i} ,
qui est un \'el\'ement de \m{\HHom(E_0,E_0\ot\ki_Y)(U_{ijk})} .
Pour obtenir une extension de $\ke$ \`a \m{X\times C_n} on peut remplacer les
\m{\Lambda_{ij}} par \
\m{\Lambda'_{ij}=\Lambda_{ij}-\beta_{ij}} ,
avec \m{\beta_{ij}} nul sur $Y$. On peut donc consid\'erer les
\m{\beta_{ij}} comme des \'el\'ements de 
\m{\HHom(\ko_C\ot\C^r,\ki_Y\ot\C^r)(U_{ij})}.
Soit \
\m{\rho'_{ijk}=\Lambda'_{jk}\Lambda'_{ij}-\Lambda'_{ik}} .
Alors on a \
\m{\rho'_{ijk}=\rho_{ijk}-\beta_{jk}\Lambda_{ij}-
\Lambda_{jk}\beta_{ij}+\beta_{ik}} .
Posons \ \m{\alpha_{ij}=(\lambda_j)^{-1}\beta_{ij}\lambda_i}, qui est un
\'el\'ement de \m{\HHom(E_0,E_0\ot\ki_Y)(U_{ij})} . Alors on a
\m{\rho'_{ijk}=0} si et seulement si
\[
(*) \ \ \ \ \ \ \mu_{ijk} \ = \ \alpha_{ij}+\alpha_{jk}-\alpha_{ik} .\]
On a \
\m{\Lambda_{kl}\rho_{ijk}-\rho_{ijl}+\rho_{ikl}-\rho_{jkl}\Lambda_{ij}=0} ,
d'o\`u il d\'ecoule que
\m{\mu_{ijk}-\mu_{ijl}+\mu_{ikl}-\mu_{jkl}=0} ,
c'est-\`a-dire que \m{(\mu_{ijk})} est un cocycle associ\'e au faisceau
\m{\HHom(E_0,E_0\ot\ki_Y)} sur \m{X\times C} et au recouvrement \m{(U_i)}. Comme
on l'a vu, on a \ \m{H^2(\HHom(E_0,E_0\ot\ki_Y))=\nsp} ,
d'o\`u l'existence des \m{\alpha_{ij}} satisfaisant l'\'egalit\'e $(*)$ et des
\m{\beta_{ij}} d\'efinissant le prolongement voulu de $\ke$.
\end{proof}

\end{sub}

\sepsub

\Ssect{Groupe de Picard}{pic_cn}

Soient \m{C_n} une courbe multiple primitive de multiplicit\'e $n$ de courbe
r\'eduite associ\'ee $C$ projective et $d$ un entier. Soit $L$ le
fibr\'e en droites sur $C$ associ\'e \`a \m{C_n}. Contrairement \`a ce qui est
affirm\'e dans \cite{dr2}, 3.2, l'existence d'un fibr\'e de Poincar\'e sur
\m{\Pic^d(C_n)\times C_n} n'est assur\'ee que sous certaines hypoth\`eses, par
exemple si \m{\deg(L)\leq 0} et si \m{L^i\not=\ko_C} pour \m{1\leq i\leq n}.
Cela est d\`u au fait qu'en g\'en\'eral les faisceaux \m{\ko_i} peuvent avoir
des sections non triviales (cf. \cite{fga}). On va d\'ecrire ici sommairement la
construction de \cite{dr2}, 3.2, qui donne en fait une vari\'et\'e lisse
\m{\Lambda_n^d(C_n)} au dessus \m{\Pic^d(C_n)} et un ``fibr\'e de Poincar\'e''
$\kd_n$ sur \m{\Lambda_n^d(C_n)\times C_n}.

Soient \m{\D_{n-1}} un fibr\'e en droites sur \m{C_{n-1}} et
\m{D=\D_{n-1\mid C}}. Si $\D$ est un fibr\'e en droites sur \m{C_n} prolongeant
\m{\D_{n-1}}, on a une suite exacte
\[0\lra D\ot L^{n-1}\lra\D\lra\D_{n-1}\lra 0 .\]
On a une suite exacte
\xmat{0\ar[r] & H^1(L^{n-1})\ar[rr]\fleq[d] & &
\Ext^1_{\ko_n}(\D_{n-1},L^{n-1}\ot
D)\ar[r]^-\phi\fleq[d] & \C\ar[r] & 0\\
 & \Ext^1_{\ko_C}(\D_{n-1\mid C},D\ot L^{n-1}) & & \Hom(\D^{(1)}_{n-1}\ot L,
D\ot L^{n-1})}
(cf. \cite{dr1} ou \ref{Ext_ii}). Les fibr\'es en droites sur \m{C_n}
prolongeant \m{\D_{n-1}} sont les extensions param\'etr\'ees par
\m{\phi^{-1}(1)}.

La vari\'et\'e \m{\Lambda_n^d(C_n)} et le fibr\'e en droites $\kd_n$ se
construisent par r\'ecurrence sur $n$ au moyen d'une ``extension universelle'':
\m{\Lambda_n^d(C_n)} est un fibr\'e en espaces affines sur \m{\Lambda_
{n-1}^d(C_{n-1})}.

\end{sub}

\sepsec

\section{Propri\'et\'es \'el\'ementaires des faisceaux coh\'erents sur les
courbes multiples primitives}\label{FCCMP}

On consid\`ere dans ce chapitre une courbe multiple primitive \m{C_n} de courbe
r\'eduite associ\'ee $C$. On utilise les notations de \ref{def_nota}.

\sepsub

\Ssect{Filtrations canoniques}{filtcan}

Soient $P$ un point ferm\'e de $C$, $M$ un \m{\ko_{nP}}-module de type fini.
Soit $\ke$ un faisceau coh\'erent sur \m{C_n}.

\sepsubsub

\begin{subsub}\label{QLL-def1}Premi\`ere filtration canonique - \rm
On d\'efinit la {\em premi\`ere filtration canonique de $M$} : c'est la
filtration
\[M_n=\nsp\subset M_{n-1}\subset\cdots\subset M_{1}\subset M_0=M\]
telle que pour \m{0\leq i< n}, \m{M_{i+1}} soit le noyau du morphisme
canonique surjectif \Nligne
\m{M_{i}\to M_{i}\ot_{\ko_{n,P}}\ko_{C,P}} .
On a donc
\[M_{i}/M_{i+1} \ = \ M_{i}\ot_{\ko_{n,P}}\ko_{C,P}, \ \ \ \
M/M_i \ \simeq \ M\ot_{\ko_{n,P}}\ko_{i,P}, \ \ \ \
M_i \ = \ z^iM .\]
 On pose, si $i>0$,
\m{G_i(M)  =  M_i/M_{i+1}} .
Le gradu\'e
\[{\rm Gr}(M) \ = \ \bigoplus_{i=0}^{n-1}G_i(M) \ = \
\bigoplus_{i=0}^{n-1}z^iM/z^{i+1}M\]
est un \m{\ko_{C,P}}-module. Les propri\'et\'es suivantes sont imm\'ediates : si
\m{1<i\leq n}\Nligne
- on a $M_i=\nsp$ si et seulement si $M$ est un
$\ko_{i,P}$-module,\Nligne
- $M_i$ est un $\ko_{n-i,P}$-module, et sa filtration canonique est \
\m{\nsp\subset M_n\subset\cdots\subset M_{i+1}\subset M_i} ,\Nligne
- tout morphisme de \m{\ko_{n,P}}-modules envoie la premi\`ere
filtration canonique du premier module sur celle du second. 

On d\'efinit de m\^eme la {\em premi\`ere filtration canonique de $\ke$} :
c'est la filtration
\[\ke_n=0\subset \ke_{n-1}\subset\cdots\subset \ke_{1}\subset \ke_0=\ke\]
telle que pour \m{0\leq i< n}, \m{\ke_{i+1}} soit le noyau du morphisme
canonique surjectif \ \m{\ke_i\to\ke_{i\mid C}}.
On a donc \
\m{\ke_{i}/\ke_{i+1}=\ke_{i\mid C}} ,
\m{\ke/\ke_i=\ke_{\mid C_i}} .
 On pose, si $i\geq 0$,
\[G_i(\ke) \ = \ \ke_i/\ke_{i+1} .\]
Le gradu\'e
\m{{\rm Gr}(\ke)} est un \m{\ko_{C}}-module. Les propri\'et\'es suivantes sont
imm\'ediates : si \m{1<i\leq n}\Nligne
- on a \m{\ke_i=\ki_C^i\ke}, et donc \ \m{{\rm
Gr}(\ke)=\bigoplus_{j=0}^{n-1}\ki_C^j\ke/\ki_C^{j+1}\ke} .\Nligne
- on a $\ke_i=0$ si et seulement si $\ke$ est un faisceau sur $C_i$,\Nligne
- $\ke_i$ est un faisceau sur $C_{n-i}$, et sa filtration canonique
est \
\m{0\subset\ke_n\subset\cdots\subset\ke_{i+1}\subset\ke_i} .\Nligne
- tout morphisme de faisceaux coh\'erents sur \m{C_n} envoie la
premi\`ere filtration canonique du premier sur celle du second.
\end{subsub}

\sepsubsub

\begin{subsub}\label{type_complet}Type complet d'un faisceau coh\'erent - \rm La
paire \
\[\Big(\big(\rg(G_0(\ke)),\ldots,\rg(G_{n-1}(\ke))\big),\big(\deg(G_0(\ke)),
\ldots , \deg(G_{n-1}(\ke)\big)\Big)\]
s'appelle le {\em type complet} de $\ke$.
\end{subsub}

\sepsubsub

\begin{subsub}\label{2-fc}Seconde filtration canonique - \rm
On d\'efinit la {\em seconde filtration canonique de $M$} : c'est la filtration
\[M^{(0)}=\nsp\subset M^{(1)}\subset\cdots\subset M^{(n-1)}\subset M^{(n)}=M\]
avec \
\m{M^{(i)} \ = \ \big\lbrace u\in M ; z^iu=0\big\rbrace} .
Si \ \m{M_n=\nsp\subset M_{n-1}\subset\cdots\subset M_1\subset M_0=M} \ est
la (premi\`ere) filtration canonique de $M$ on a \ \m{M_i\subset M^{(n-i)}} \
pour \m{0\leq i\leq n}. On pose, si $i>0$,
\m{G^{(i)}(M)  =  M^{(i)}/M^{(i-1)}} .
Le gradu\'e
\[{\rm Gr}_2(M) \ = \ \bigoplus_{i=1}^nG^{(i)}(M) \]
est un \m{\ko_{C,P}}-module.
Les propri\'et\'es suivantes sont imm\'ediates : si \m{1<i\leq n}\Nligne
- $M^{(i)}$ est un $\ko_{i,P}$-module, et sa filtration canonique
est \hfil\break
\m{\nsp\subset M^{(1)}\subset\cdots\subset M^{(i-1)}\subset M^{(i)}}~,
\Nligne
- tout morphisme de \m{\ko_{n,P}}-modules envoie la seconde
filtration canonique du premier module sur celle du second. 

On d\'efinit de m\^eme la {\em seconde filtration canonique de $\ke$} :
\[\ke^{(0)}=\nsp\subset \ke^{(1)}\subset\cdots\subset
\ke^{(n-1)}\subset \ke^{(n)}=\ke  .\]
 On pose, si $i>0$,
\[G^{(i)}(\ke) \ = \ \ke^{(i)}/\ke^{(i-1)} .\]
Le gradu\'e
\m{{\rm Gr}_2(\ke)} est un \m{\ko_{C}}-module. Les propri\'et\'es suivantes sont
imm\'ediates : si \m{0< i\leq n}\Nligne
- $\ke^{(i)}$ est un faisceau sur $C_i$, et sa filtration
canonique est \
\m{0\subset\ke^{(1)}\subset\cdots\subset\ke^{(i-1)}\subset\ke^{(i)}} ,
\Nligne
- tout morphisme de faisceaux coh\'erents sur \m{C_n} envoie la
seconde filtration canonique du premier sur celle du second.
\end{subsub}

\sepprop

\begin{subsub}\label{filt1}{\bf Proposition : } On a des isomorphismes
canoniques :

(i) \m{\ke_i\simeq(\ke/\ke^{(i)})\ot\L^i} .

(ii) \m{(\ke^{(i)})^{(j)}\simeq\ke^{(\min(i,j))}} .

(iii) \m{(\ke_i)_j\simeq\ke_{i+j}} .

(iv) \m{(\ke_i)^{(j)}\simeq(\ke^{(i+j)}/\ke^{(i)})\ot\L^i} .

(v) \m{(\ke^{(i)})_j=0} \ si \m{i\leq j}, et \
\m{(\ke^{(i)})_j\simeq(\ke_j)^{(i-j)}} \ si \m{i>j}.
\end{subsub}

\begin{proof} Imm\'ediat, en examinant ce qui se passe en chaque point de $C$.
\end{proof}

\sepprop

\begin{subsub}\label{rel_filt}Relations entre les deux filtrations - \rm
Soient $i$, $k$ des entiers tels que \m{0<i<n}, \m{0< k\leq i+1}.
Le morphisme canonique \m{\ke\ot\L\to\ke} induit un morphisme de faisceaux
coh\'erents sur $C$ :
\[\lambda_{ik}=\lambda_{ik}(\ke):G^{(i+1)}(\ke)\ot L^k\lra G^{(i+1-k)}(\ke) .\]
On posera, si \m{0\leq i<n},
\[\Gamma^{(i)}(\ke) \ = \ \coker(\lambda_{i+1,1}) .\]
Soient $j$, $m$ des entiers tels que \m{j,m>0}, \m{j+m\leq n}. Le morphisme
canonique \m{\ke\ot\L\to\ke} induit un morphisme de faisceaux coh\'erents sur
$C$ :
\[\mu_{jm}=\mu_{jm}(\ke):G_j(\ke)\ot L^m\lra G_{j+m}(\ke) .\]
On posera, si \m{0\leq j<n},
\[\Gamma_j(\ke) \ = \ \ker(\mu_{j1}) .\]
\end{subsub}

\sepprop

\begin{subsub}\label{filt2_l}{\bf Lemme : } Le noyau du morphisme canonique
surjectif \ \m{\ke_i\ot\L\to\ke_{i+1}} \ est isomorphe \`a \m{G^{(i+1)}(\ke)\ot
L^{i+1}} .

\end{subsub}

\begin{proof} Cela d\'ecoule du diagramme commutatif avec lignes exactes
\xmat{0\ar[r] & \ke^{(i)}\ar[r]\flinc[d] & \ke\ar[r]\fleq[d] &
\ke_i\ot\L^{-i}\ar[r]\flon[d] & 0\\
0\ar[r] & \ke^{(i+1)}\ar[r] & \ke\ar[r] & \ke_{i+1}\ot\L^{-i-1}\ar[r] & 0
}
\end{proof}

\sepprop

\begin{subsub}\label{filt2}{\bf Proposition : } 1 - Le morphisme de faisceaux
\m{\lambda_{ik}} est injectif. De plus, on a une suite exacte canonique
\[0\lra L^{i+1-k}\ot\coker(\lambda_{i,k})\lra\ke_{i-k\mid C_{k}}\ot
\L\lra\ke_{i+1-k}/\ke_{i+1}\lra 0 .\]
2 - Le morphisme \m{\mu_{jm}} est surjectif, et on a une suite exacte de
faisceaux coh\'erents sur \m{C_m} :
\[0\lra(\ke^{(j+m+1)}/\ke^{(j+1)})\ot\L^{j+m+1}\lra(\ke^{(j+m)}/\ke^{(j)})\ot
\L^{j+m}\lra\ker(\mu_{jm})\lra 0 .\]
\end{subsub}

\begin{proof} On ne d\'emontrera que 1-, 2- \'etant analogue.
L'injectivit\'e de \m{\lambda_{ik}} se d\'emontre ais\'ement en examinant ce
qui se passe en chaque point. La suite exacte est obtenue \`a partir du
diagramme commutatif avec lignes exactes d\'eduites du lemme \ref{filt2_l}
\xmat{
0\ar[r] & G^{(i+1)}(\ke)\ot L^{i+1}\ar[r]\flinc[d]^{\lambda_{ik}\ot
I_{L^{i+1-k}}} & \ke_i\ot\L\ar[r]\flinc[d] & \ke_{i+1}\ar[r]\flinc[d] & 0 \\
0\ar[r] & G^{(i+1-k)}(\ke)\ot L^{i+1-k}\ar[r] & \ke_{i-k}\ot\L\ar[r]
& \ke_{i+1-k}\ar[r] & 0
}
\end{proof}

\sepprop

\begin{subsub}\label{filt3}{\bf Corollaire : } Si \ \m{0\leq i<n}, on a un
isomorphisme canonique \ \m{\Gamma_i(\ke)\simeq\Gamma^{(i)}(\ke)\ot L^{i+1}}.
\end{subsub}

\begin{proof} Cela d\'ecoule de la proposition \ref{filt2}, 1- en prenant
\m{k=1}, ou de 2- en prenant \m{m=1}, ou directement du lemme \ref{filt2_l}.
\end{proof}

\end{sub}

\sepprop

\begin{subsub}\label{inj_surj}{\bf Proposition : } Soient $\ke$, $\kf$ des
faisceaux coh\'erents sur \m{C_n}. Soit \ \m{\phi:\ke\to\kf} \ un morphisme. 
Alors $\phi$ est surjectif si et seulement si le morphisme \ \m{\phi_{\mid
C}:G_0(\ke)\to G_0(\kf)} \ de faisceaux sur $C$ induit par $\phi$ l'est. Si
c'est le cas les morphismes \ \m{G_i(\ke)\to G_i(\kf)} induits par $\phi$ sont
surjectifs pour \m{1\leq i<n}.

Le morphisme $\phi$ est injectif si et seulement si le morphisme \
\m{G^{(1)}(\ke)=\ke^{(1)}\to G^{(1)}(\kf)=\kf^{(1)}} \ de faisceaux sur $C$
l'est. Si c'est le cas les morphismes \ \m{G^{(i)}(\ke)\to G^{(i)}(\kf)} \
induits par $\phi$ sont injectifs pour \m{1\leq i<n}.
\end{subsub}

\begin{proof}
Supposons que $\phi$ soit surjectif. Alors \m{\phi_{\mid C}} l'est. La
r\'eciproque se d\'emontre par r\'ecurrence sur $n$, le cas \m{n=1} \'etant
trivial. Supposons qu'elle soit vraie pour \m{n-1\geq 1}, et que \m{\phi_{\mid
C}} soit surjectif. Alors $\phi$ induit \ \m{\phi_1:\ke_1\to\kf_1} . On a un
diagramme commutatif
\xmat{
G_0(\ke)\ot L\flon[rr]^-{\mu_{01}(\ke)}\flon[d] & & G_1(\ke)=G_0(\ke_1)\ar[d]\\
G_0(\kf)\ot L\flon[rr]^-{\mu_{01}(\kf)} & & G_1(\kf)=G_0(\kf_1)}
d'o\`u il d\'ecoule que \m{\phi_{1\mid C}} est surjectif. D'apr\`es
l'hypoth\`ese de r\'ecurrence, $\phi_1$ est surjectif. La surjectivit\'e de
$\phi$ se d\'eduit ais\'ement de celles de \m{\phi_{\mid C}} et \m{\phi_1}.
La surjectivit\'e du morphisme \ \m{G_i(\ke)\to G_i(\kf)} \ d\'ecoule du
diagramme commutatif
\xmat{
G_0(\ke)\ot L^i\flon[rr]^-{\mu_{0i}(\ke)}\flon[d] & & G_i(\ke)\ar[d]\\
G_0(\kf)\ot L^i\flon[rr]^-{\mu_{0i}(\kf)} & & G_i(\kf)}

La seconde assertion se d\'emontre de mani\`ere analogue.
\end{proof}

\sepsub

\Ssect{Invariants}{inva}

\begin{subsub}Rang g\'en\'eralis\'e - \rm L'entier \
\m{R(M)=\rg({\rm Gr}(M))} \
s'appelle le {\em rang g\'en\'eralis\'e} de $M$.

L'entier \
\m{R(\ke)=\rg({\rm Gr}(\ke))} \
s'appelle le {\em rang g\'en\'eralis\'e} de $\ke$. On a donc
\m{R(\ke)=R(\ke_P)} pour tout \m{P\in C}.

\end{subsub}

\sepsubsub

\begin{subsub}Degr\'e g\'en\'eralis\'e - \rm L'entier \
\m{\Deg(\ke)=\deg({\rm Gr}(\ke))} \
s'appelle le {\em degr\'e g\'en\'eralis\'e de } $\ke$.

Si \m{R(\ke)>0} on pose \m{\mu(\ke)=\Deg(\ke)/R(\ke)} et on appelle ce nombre la
{\em pente} de $\ke$.
\end{subsub}

\sepprop

\begin{subsub}{\bf Proposition : }\label{str_gen2}
1 - Soit \
\m{0\to M'\to M\to M''\to 0} \
une suite exacte de \m{\ko_{n,P}}-modules de type fini. Alors on a \
\m{R(M)=R(M')+R(M'')} .

2 - Soit \
\m{0\to E\to F\to G\to 0} \
une suite exacte de faisceaux coh\'erents sur \m{C_n}. Alors on a \
\m{R(F)=R(E)+R(G)} , \m{\Deg(F)=\Deg(E)+\Deg(G)} .
\end{subsub}

\sepprop

\begin{subsub}\label{filt4}{\bf Proposition :} Posons \ \m{r_j=\rg(G_j(\ke))} \
pour \m{0\leq j<n}. On a, pour \m{0\leq i<n},
\[\rg(G^{(i+1)}(\ke))= \ \rg(G_i(\ke))\]
et
\[\deg(G^{(i+1)}(\ke)) \ = \
\deg(G_i(\ke))+(r_{i+1}+\cdots+r_{n-1}-i.r_i)\deg(L) .\]

\end{subsub}

\begin{proof} 
Posons, pour \m{0\leq j<n}, \m{s_j=\rg(G^{(j+1)}(\ke))}, \m{d_j=\rg(G_j(\ke))},
\hfil\break\m{e_j=\rg(G^{(j+1)}(\ke))}, et si \m{j<n-1},
\m{\delta_j=d_{j+1}-d_j},
\m{\epsilon_j=e_{j+1}-e_j}. On a, d'apr\`es le corollaire \ref{filt3}, une suite
exacte
\[0\lra G^{(j+2)}(\ke)\ot L^{j+2}\lra G^{(j+1)}(\ke)\ot L^{j+1}\lra G_j(\ke)\ot
L\lra G_{j+1}(\ke)\lra 0 .\]
On en d\'eduit que  \ \m{s_j-s_{j+1}=r_j-r_{j+1}}, et on en d\'eduit la
premi\`ere \'egalit\'e de la proposition, compte tenu du fait que
\[r_0+\cdots+\cdots r_{n-1} \ = \ s_0+\cdots+\cdots s_{n-1} \ = \ R(\ke)\]
(d'apr\`es la proposition \ref{str_gen2}). On a, d'apr\`es la suite exacte
pr\'ec\'edente
\begin{equation}\label{eq2}
\epsilon_j-\delta_j \ = \ \big(j.r_j-(j+2)r_{j+1}\big)\deg(L) .
\end{equation}
On a \ \m{d_{i+1}=d_0+\delta_0+\cdots\delta_i},
\m{e_{i+1}=e_0+\epsilon_0+\cdots\epsilon_i}, et 
\begin{eqnarray*}
 \sigg_{i=0}^{n-1}d_i & = & nd_0+(n-1)\delta_0+(n-2)\delta_1+\cdots+
\delta_{n-2}\\
& = & \Deg(\ke)  =  \sigg_{i=0}^{n-1}e_i\\
& = & ne_0+(n-1)\epsilon_0+(n-2)\epsilon_1+\cdots+\epsilon_{n-2} ,
\end{eqnarray*}
d'o\`u, en utilisant (\ref{eq2})
\begin{eqnarray*}
n(d_0-e_0) & = &
\big(((n-1)(0.r_0-2r_1)+\cdots+(n-j)\big((j-1)r_{j-1}-(j+1)r_j\big)+\cdots
\big)\deg(L)\\
& = & -n(r_1+\cdots+r_{n-1})\deg(L) .
\end{eqnarray*}
Donc \ \m{e_0=d_0+(r_1+\cdots+r_{n-1})\deg(L)}, c'est la seconde
\'egalit\'e de la proposition \ref{filt4} pour \m{i=0}. Le cas g\'en\'eral
se d\'emontre ensuite par r\'ecurrence sur $i$ en utilisant (\ref{eq2}).
\end{proof}
\end{sub}

\sepsub

\Ssect{Torsion}{tors}

Soit $M$ un \m{\ko_{nP}}-module de type fini. Le {\em sous-module de torsion
\m{T(M)} de $M$} est constitu\'e des \'el\'ements annul\'es par une puissance de
$x$. On dit que $M$ est {\em sans torsion} si ce sous-module est nul. C'est donc
le cas si et seulement si pour tout \m{m\in M} non nul et tout entier \m{p>0} on
a \m{x^pm\not=0}.

Soit $\ke$ un faisceau coh\'erent sur \m{C_n}. Le {\em sous-faisceau de torsion
\m{T(\ke)} de $\ke$} est le sous-faisceau maximal de $\ke$ dont le support est
fini. Pour tout point ferm\'e $P$ de $C$ on a \ \m{T(\ke)_P=T(\ke_P)}. On dit
que $\ke$ est un faisceau {\em de torsion} si \ \m{\ke=T(\ke)} , ou ce qui
revient au m\^eme, si son support est fini.

\sepprop

\begin{subsub}\label{prop_tor}{\bf Proposition : } Les conditions suivantes
sont \'equivalentes :
\begin{enumerate}
\item[(i)] $\ke$ est sans torsion.
\item[(ii)] $G^{(1)}(\ke)=\ke^{(1)}$ est localement libre.
\end{enumerate}
De plus, si $\ke$ est sans torsion, tous les faisceaux $G^{(i)}(\ke)$ sur $C$
sont localement libres.
\end{subsub}

\begin{proof}
Supposons que $\ke$ soit sans torsion. Il en est de m\^eme de tous ses
sous-faisceaux, et donc de \m{G^{(1)}(\ke)}. Comme les faisceaux sans torsion
sur $C$ sont les faisceaux localement libres, \m{G^{(1)}(\ke)} est localement
libre.

R\'eciproquement, supposons que \m{T(\ke)\not=0}. Il existe alors un point
ferm\'e $P$ de $C$ et \m{u\in\ke_P} non nul tel qu'il existe un entier \m{p>0}
tel que \m{x^pu=0}. Soit $q$ le plus grand entier tel que \m{z^qu\not=0}. Alors
on a \m{z^qu\in G^{(1)}(\ke_P)}, et \m{x^p.z^qu=0}, donc \m{G^{(1)}(\ke)_P} a
de la torsion.

Le d\'emonstration de la seconde assertion est analogue.
\end{proof}

\end{sub}

\sepsub

\Ssect{Faisceaux coh\'erents quasi localement libres}{QLL}

Soit $M$ un \m{\ko_{n,P}}-module de type fini.
On dit que $M$ est {\em quasi libre} s'il existe des entiers \m{m_1,\ldots,m_n}
non n\'egatifs et un isomorphisme
\m{M\simeq\bigoplus_{i=1}^n m_i\ko_{i,P}} .
Les entiers \m{m_1,\ldots,m_n} sont uniquement d\'etermin\'es.
On dit que $M$ est {\em de type} \m{(m_1,\ldots,m_n)}. On a \
\m{R(M)=\sigg_{i=1}^ni.m_i}  .

Soit $\ke$ un faisceau coh\'erent sur \m{C_n}.
On dit que $\ke$ est {\em quasi localement libre} en un point $P$ de
$C$ s'il existe un ouvert $U$ de \m{C_n} contenant
$P$ et des entiers non n\'egatifs \m{m_1,\ldots,m_n} tels que pour tout
point $Q$ de $U$, \m{\ke_{n,Q}} soit quasi localement libre de type
\m{m_1,\ldots,m_n}. Les entiers \m{m_1,\ldots,m_n} sont uniquement
d\'etermin\'es et ne d\'ependent que de $\ke$,
et on dit que \m{(m_1,\ldots,m_n)} est le {\em type de } $\ke$.

On dit que $\ke$ est {\em quasi localement libre} s'il l'est en tout point de
\m{C_n}.

\sepprop

\begin{subsub}{\bf Th\'eor\`eme : }\label{str_gen0}
Soient $\ke$ un  faisceau coh\'erent sur \m{C_n} et $P$ un point ferm\'e de $C$.
Alors les conditions suivantes sont \'equivalentes : 

(i) $\ke_P$ est quasi localement libre en $P$.

(ii) Pour \m{0\leq i<n}, \m{G_i(\ke)} est libre en $P$.

(iii) Pour \m{0\leq i<n}, \m{\Gamma^{(i)}(\ke)} est libre en $P$.
\end{subsub}

\begin{proof}
L'\'equivalence de (i) et (ii) est d\'emontr\'ee dans \cite{dr2}, th\'eor\`eme
5.1.4. Montrons que (ii) et (iii) sont \'equivalentes : d'apr\`es \ref{filt3},
(iii) \'equivaut \`a

(iii)' Pour \m{0\leq i<n}, \m{\Gamma_i(\ke)} est libre en $P$.

L'\'equivalence de (ii) et (iii)' est imm\'ediate, compte tenu de la
surjectivit\'e des morphismes \m{\mu_{i1}} et du fait que \
\m{\Gamma_{n-1}(\ke)=G_{n-1}(\ke)\ot L} .
\end{proof}

\sepprop

\begin{subsub}\label{coro_str}{\bf Corollaire : } Soit $\ke$ un  faisceau
coh\'erent sur \m{C_n}. Alors les conditions suivantes sont \'equivalentes : 

(i) $\ke$ est quasi localement libre.

(ii) Pour \m{0\leq i<n}, \m{G_i(\ke)} est localement libre sur $C$.

(iii) Pour \m{0\leq i<n}, \m{\Gamma^{(i)}(\ke)} est localement libre sur $C$.
\end{subsub}

\sepprop

\begin{subsub}\label{coro_str2}{\bf Corollaire : } Soit $\ke$ un  faisceau
coh\'erent sans torsion sur \m{C_n}. Alors les conditions suivantes sont
\'equivalentes : 

(i) $\ke$ est localement libre.

(ii) $\ke$ admet une r\'esolution localement libre de longueur finie.
\end{subsub}

\begin{proof}
Supposons que $\ke$ admette une r\'esolution localement libre de longueur finie.
Alors on a \ \m{\Tor^k(\ke,\kf)=\nsp} \ pour tout $k$ suffisamment grand et tout
faisceau coh\'erent $\kf$ sur \m{C_n}. Soit $i$ un entier tel que \m{1\leq i<n}.
On consid\`ere la r\'esolution localement libre \'evidente de $\ko_i$
\xmat{\cdots\L^{2n}\ar[r] & \L^{n+i}\ar[r] & \L^n\ar[r] & \L^i\ar[r] &
\ko_n\ar[r] & \ko_i ,}
d\'efinie par les morphismes canoniques \m{\L^i\to\ko_n} et
\m{\L^{n-i}\to\ko_n}. La suite induite
\xmat{\cdots\L^{2n}\ot\ke\ar[r] & \L^{n+i}\ot\ke\ar[r] & \L^n\ot\ke\ar[r] &
\L^i\ot\ke\ar[r] & \ko_n\ot\ke}
a pour cohomologies les \m{\Tor^k(\ke,\ko_i)}. Elle est donc exacte \`a partir
d'un certain degr\'e. Il en d\'ecoule que \ \m{\ke^{(i)}=\ke_{n-i}} . Puisque
$\ke$ est sans torsion, les \m{G^{(i)}(\ke)} sont localement libres d'apr\`es la
proposition \ref{prop_tor}. Il en est donc de m\^eme des \m{G_i(\ke)}. Donc
d'apr\`es le th\'eor\`eme \ref{str_gen0}, $\ke$ est quasi localement libre.
Soient \m{(m_1,\ldots,m_n)} le type de $\ke$. Alors en un point $P$ de $C$ o\`u
$\ke_P$ est identifi\'e \`a \m{\OPL_{i=1}^nm_i\ko_{iP}} on a
\[\OPL_{i=1}^{n-1}m_i\ko_{iP} \ \subset \ \ke_P^{(n-1)} .\]
Mais si \m{m_1,\ldots,m_{n-1}} ne sont pas tous nuls on a
\[\OPL_{i=1}^{n-1}m_i\ko_{iP} \ \not\subset \ \ke_{1P} .\]
On a donc \ \m{m_1=\cdots=m_{n-1}=0} , et $\ke$ est localement libre.
\end{proof}

\sepprop

Le r\'esultat suivant sera utilis\'e au chapitre \ref{fam_qll} :

\sepprop

\begin{subsub}\label{f_end}{\bf Proposition : } Soit $\ke$ un faisceau quasi
localement libre sur \m{C_n}. On note \ \m{r_j=\rg(G_j(\ke))} \ pour \m{0\leq
j<n}. Alors le faisceau des endomorphismes \m{\EEnd(\ke)} est aussi quasi
localement libre, et on a
\[R(\EEnd(\ke))=\sigg_{i=0}^{n-1}r_i^2 \ , \quad\quad
\Deg(\EEnd(\ke))=\big (\sigg_{0\leq i<j<n}r_ir_j\big)\deg(L) \ .\]
\end{subsub}

\begin{proof} Le fait que \m{\EEnd(\ke)} est quasi localement libre est
imm\'ediat. On a une suite exacte canonique
\[0\lra\HHom(\ke_{\mid C},\ke^{(1)})\lra\EEnd(\ke)\lra\EEnd(\ke_1)\lra 0 ,\]
d'o\`u on d\'eduit que
\[R(\EEnd(\ke))=R(\EEnd(\ke_1))+\rg((\ke_{\mid C})^*\ot\ke^{(1)}) \ ,
\Deg(\EEnd(\ke))=\Deg(\EEnd(\ke_1))+\deg((\ke_{\mid C})^*\ot\ke^{(1)}) \ .\]
En raisonnant par r\'ecurrence sur $n$ et en utilisant la proposition
\ref{filt1}, (iv), on obtient les formules
\[R(\EEnd(\ke))=\sigg_{0<i<n}\rg(G_i(\ke)^*\ot G^{(i+1)}(\ke)) \ , \
\Deg(\EEnd(\ke))=\sigg_{0<i<n}\deg(G_i(\ke)^*\ot G^{(i+1)}(\ke)\ot L^i) \ .\]
Le r\'esultat d\'ecoule alors ais\'ement de la proposition \ref{filt4}.
\end{proof}

\end{sub}

\sepsub

\Ssect{Construction des faisceaux coh\'erents}{constr}

On d\'ecrit ici le moyen de construire un faisceau coh\'erent $\ke$ sur
\m{C_n}, connaissant \m{\ke_{\mid C}} et \m{\ke_1} (ou \m{\ke^{(1)}} et
\m{\ke_1}), qui sont des faisceaux sur $C$ et \m{C_{n-1}} respectivement.
Cela de permet de faire des d\'emonstrations par r\'ecurrence sur $n$ dans les
chapitres suivants.

On consid\`ere dans ce chapitre une courbe multiple primitive \m{C_n} de courbe
r\'eduite associ\'ee $C$. On utilise les notations de \ref{def_nota}. Si $\kg$
est un fibr\'e vectoriel sur \m{C_k} (\m{1\leq k<n}), il existe un fibr\'e
vectoriel $\G$ sur \m{C_n} qui prolonge $\kg$ (cf. \cite{dr2}, th\'eor\`eme
3.1.1). En particulier le faisceau d'id\'eaux \m{\ki_C} de $C$ dans \m{C_n} se
prolonge en un fibr\'e en droites $\L$ sur \m{C_n}. On en d\'eduit une
r\'esolution localement libre de $\kg$ :
\[\cdots\G\ot\L^{n+k}\lra\G\ot\L^n\lra\G\ot\L^k\lra\G\lra\kg\lra 0 .\]

\sepsub

\begin{subsub}\label{Ext_i} Premi\`ere construction - \rm Soient $\kf$ un
faisceau coh\'erent sur \m{C_{n-1}} et $E$ un fibr\'e vectoriel sur $C$. On
s'int\'eresse aux faisceaux coh\'erents $\ke$ sur \m{C_n} tels que
\m{\ke_{\mid C}=E} et \m{\ke_1=\kf}. Soit \
\m{0\lra\kf\lra\ke\lra E\lra 0} \ une suite exacte, associ\'ee \`a \
\m{\sigma\in\Ext_{\ko_n}^1(E,\kf)}.
Soit \ \m{\pi_\ke:\ke\ot\ki_C\to\ke} \ le morphisme canonique. On a
\[\imm(\pi_\ke)\subset\kf , \quad \pi_\ke(\kf)=\imm(\pi_\kf) .\]
Donc \m{\pi_\ke} induit un morphisme
\[\Phi_{\kf,E}(\sigma) :E\ot L\lra\kf_{\mid C} ,\]
et on obtient ais\'ement le
\end{subsub}

\sepprop

\begin{subsub}\label{le_i_1}{\bf Lemme : } On a \ \m{\ke_{\mid C}=E} \
et \ \m{\ke_1=\kf} \ si et seulement si \m{\Phi_{\kf,E}(\sigma)} est
surjectif.
\end{subsub}

\sepprop

\begin{subsub}\label{pro_i_1}{\bf Proposition : } On a une suite exacte
canonique
\xmat{0\ar[r] & \Ext^1_{\ko_C}(E,\kf^{(1)})\ar[r] & \Ext^1_{\ko_n}(E,\kf)
\ar[rr]^-{\Phi_{\kf E}} & & \Hom(E\ot L,\kf_{\mid C})\ar[r] & 0 .}
\end{subsub}

\begin{proof}
Soit $\E$ un fibr\'e vectoriel sur \m{C_n} prolongeant $E$ et
\[\cdots\E\ot\L^{n+1}\lra\E\ot\L^n\lra\E\ot\L\lra\E\lra E\lra 0\]
la r\'esolution localement libre induite. On en d\'eduit l'isomorphisme
canonique
\[\EExt^1_{\ko_n}(E,\kf) \ \simeq \ \HHom(E\ot L,\kf_{\mid C}) .\]
Le r\'esultat d\'ecoule alors de la suite exacte canonique
\[0\lra H^1(\HHom(E,\kf))\lra\Ext^1_{\ko_n}(E,\kf)\lra
H^0(\EExt^1_{\ko_n}(E,\kf))\lra 0 .\]
d\'eduite de la suite spectrale des Ext (cf. \cite{go}, 7.3).
\end{proof}

\sepsubsub

\begin{subsub}\label{Ext_ii} Seconde construction - \rm Soit $\ke$ un faisceau
quasi localement libre sur \m{C_n}. On a \ \m{\EExt^1_{\ko_n}(\ke,V)\simeq
\HHom((\ke^{(1)}/\ke_{n-1})\ot L ,V)} , d'o\`u la suite exacte
\xmat{0\ar[r] & \Ext^1_{\ko_C}(\ke_{\mid C},V)\ar[r] & 
\Ext^1_{\ko_n}(\ke,V)\ar[rr]^-{\Theta_{\ke,V}} & &
\Hom((\ke^{(1)}/\ke_{n-1})\ot L ,V)\ar[r] & 0 .}

\end{subsub}

\sepprop

\begin{subsub}\label{pro_ii_2}{\bf Proposition : } Soient $\ke$ un faisceau
coh\'erent sur \m{C_n} et $V$ un fibr\'e vectoriel sur $C$. Soient \ \m{0\to
V\to\kf\to\ke\to 0} \ une extension et \m{\sigma\in
\Ext^1_{\ko_n}(\ke,V)} \ l'\'el\'ement associ\'e.

1- Le morphisme canonique \ \m{\kf\ot\L\to\kf} \ en induit un,
\m{(\ke^{(1)}/\ke_{n-1})\ot L\to V} , qui n'est autre que
\m{\Theta_{\ke,V}(\sigma)}.

2 - Si \m{\ke} est concentr\'e sur \m{C_{n-1}}, on a \ \m{V=\kf^{(1)}} si et
seulement si \m{\Theta_{\ke,V}(\sigma)} est injectif.
\end{subsub}

\end{sub}

\sepsec

\section{Faisceaux r\'eflexifs}

On consid\`ere dans ce chapitre une courbe multiple primitive \m{C_n} de courbe
r\'eduite associ\'ee $C$. On utilise les notations de \ref{def_nota}.

\sepsub

\Ssect{Dualit\'e}{dual}

Soient \m{P\in C} et $M$ un  \m{\ko_{n,P}}-module de type fini. On note
\m{M^{\vee_n}} le {\em dual} de $M$ :\Nligne
\m{M^{\vee_n}=\Hom(M,\ko_{n,P})}~. Si aucune confusion n'est \`a craindre on
notera \m{M^\vee=M^{\vee_n}}.
Si $N$ est un \m{\ko_{C,P}}-module, on note \m{N^*} le dual de $N$ :
\m{N^*=\Hom(N,\ko_{C,P})}.

Soit $\ke$ un faisceau coh\'erent sur \m{C_n}. On note \m{\ke^{\vee_n}} le {\em
dual} de $\ke$ : \
\m{\ke^{\vee_n}=\HHom(\ke,\ko_n)} . Si aucune confusion n'est \`a craindre on
notera \m{\ke^\vee=\ke^{\vee_n}}.
Si $E$ est un faisceau coh\'erent sur $C$, on note \m{E^*} le dual de $E$ :
\m{E^*=\HHom(E,\ko_C)}.
Ces notations sont justifi\'ees par le fait que \m{E^\vee\not=E^*}, et plus
g\'en\'eralement on a

\sepprop

\begin{subsub}\label{lem1}{\bf Lemme : }Soient $i$ un entier tel que \m{1\leq
i\leq n} et $\ke$ un faisceau coh\'erent sur \m{C_i}. Alors on a un isomorphisme
canonique
\[\ke^{\vee_n} \ \simeq \ \ke^{\vee_i}\otimes\ki_C^{n-i} .\]
En particulier, pour tout faisceau coh\'erent $E$ sur $C$, on a \
\m{E^{\vee_n}\simeq E^*\ot L^{n-1}} .
\end{subsub}

\begin{proof} Imm\'ediat. \end{proof}

\sepprop

\begin{subsub}\label{x1}{\bf Proposition : } Soit $\ke$ un faisceau coh\'erent
sur \m{C_n}. Alors on a, pour \m{1\leq i<n}, \m{(\ke^\vee)^{(i)}=
(\ke/\ke_i)^\vee} .
\end{subsub}

\begin{proof} Imm\'ediat. \end{proof}

\sepprop

\begin{subsub}\label{x1b}Dualit\'e des faisceaux quasi localement libres - \rm
Soit $\ke$ un faisceau quasi localement libre sur \m{C_n}. Alors on montre
ais\'ement qu'on a, pour \m{0\leq i<n}, des isomorphismes canoniques
\[(\ke^\vee)_i\simeq(\ke_i)^\vee\ot\L^i ,\quad G^{(i+1)}(\ke^\vee)\simeq
G_i(\ke)^*\ot L^{n-1} .\]
\end{subsub}

\end{sub}

\sepsub

\Ssect{Faisceaux r\'eflexifs - caract\'erisation}{refl}

Soit $\ke$ un faisceau coh\'erent sur \m{C_n}. Rappelons qu'on dit que $\ke$ est
{\em r\'eflexif} si le morphisme canonique \ \m{\ke\to\ke^{\vee\vee}} \ est un
isomorphisme.

\sepprop

\begin{subsub}\label{lem2}{\bf Lemme : } Soit $\ke$ un faisceau quasi localement
libre sur \m{C_n}. Alors $\ke$ est r\'eflexif et on a \
\m{\EExt^i(\ke,\ko_n)=\nsp} \ pour tout \m{i\geq 1}.
\end{subsub}

\begin{proof} Il suffit de montrer que pour tout point ferm\'e $P$ de $C$ et
tout entier $m$ tel que \m{1\leq m\leq n}, \m{\ko_{mP}} est un
\m{\ko_{nP}}-module r\'eflexif et qu'on a \ \m{\Ext^i_{\ko_{nP}}(\ko_{mP},
\ko_{nP})=\nsp} \ pour tout \m{i\geq 1}. La premi\`ere assertion est imm\'ediate
et la seconde se d\'emontre en utilisant la r\'esolution libre de
\m{\ko_{mP}} :
\xmat{\cdots\ar[r] & \ko_{nP}\ar[r]^{\times z^m} & \ko_{nP}\ar[r]^{\times
z^{n-m}} & \ko_{nP}\ar[r]^{\times z^m} & \ko_{nP}\ar[r] & \ko_{mP}
}
\end{proof}

\sepprop

\begin{subsub}\label{theo_refl}{\bf Th\'eor\`eme : } Soit $\ke$ un faisceau
coh\'erent sur \m{C_n}. Alors les propri\'et\'es suivantes sont \'equivalentes :
\begin{itemize}
\item[(i)] $\ke$ est r\'eflexif.
\item[(ii)] $\ke$ est sans torsion.
\item[(iii)] On a \ $\EExt^1_{\ko_n}(\ke,\ko_n)=0$ .
\end{itemize}
Si les conditions pr\'ec\'edentes sont r\'ealis\'ees on a de plus \
\m{\EExt^i_{\ko_n}(\ke,\ko_n)=0} \ pour tout \m{i\geq 1}.
\end{subsub}

\begin{proof} On fait une r\'ecurrence sur $n$. Le r\'esultat est bien connu
pour \m{n=1}. Supposons que \m{n>1} et que le r\'esultat soit vrai pour \m{n-1}.

Soit $\kf$ un faisceau sans torsion sur \m{C_{n-1}}. D'apr\`es le lemme
\ref{lem1}, $\kf$ est aussi r\'eflexif en tant que faisceau sur \m{C_n}. On va
montrer qu'on a \ \m{\EExt^i_{\ko_n}(\kf,\ko_n)=0} \ pour \m{i\geq 1}. On
peut trouver une suite exacte de faisceaux coh\'erents sans torsion sur
\m{C_{n-1}}
\[0\lra\kn\lra\ku\lra\kf\lra 0\]
telle que $\ku$ soit quasi localement libre (par exemple on peut utiliser un
fibr\'e en droites tr\`es ample \m{\ko(1)} sur \m{C_{n-1}} et prendre pour $\ku$
une somme directe de \m{\ko(-m)} pour un $m$ suffisamment grand). On en d\'eduit
la suite exacte
\[0\lra\kf^{\vee_{n-1}}\lra\ku^{\vee_{n-1}}\lra\kn^{\vee_{n-1}}\lra 0 ,\]
et donc aussi la suivante :
\[0\lra\kf^{\vee_n}\lra\ku^{\vee_n}\lra\kn^{\vee_n}\lra 0 .\]
On en d\'eduit que
\[\EExt^1_{\ko_n}(\kf,\ko_n) \ \subset \ \EExt^1_{\ko_n}(\ku,\ko_n) .\]
D'apr\`es le lemme \ref{lem2}, on a \ \m{\EExt^1_{\ko_n}(\ku,\ko_n)=0}, donc
on a bien \ \m{\EExt^1_{\ko_n}(\kf,\ko_n)=0}. On a aussi
\[\EExt^i_{\ko_n}(\kn,\ko_n)\simeq\EExt^{i+1}_{\ko_n}(\kf,\ko_n)\]
pour \m{i\geq 1}, d'o\`u il d\'ecoule par r\'ecurrence sur $i$ que \
\m{\EExt^i_{\ko_n}(\kf,\ko_n)=0} \ pour \m{i\geq 1}.

Il est \'evident qu'un faisceau ayant de la torsion n'est pas r\'eflexif, donc
(i) entraine (ii).

Soit $\ke$ un faisceau coh\'erent sans torsion sur \m{C_n}.
Montrons que $\ke$ est r\'eflexif et que \ \m{\EExt^1_{\ko_n}(\ke,\ko_n)=0}.
On consid\`ere la suite exacte
\begin{equation}\label{eq1}
0\lra\ke^{(1)}\lra\ke\lra\ke/\ke^{(1)}=\ke_1\ot\L^\vee\lra 0 .
\end{equation}
Les faisceaux \m{\ke^{(1)}} et \m{\ke_1\ot\L^\vee} sont de support contenu dans
\m{C_{n-1}} et sans torsion. Puisque \ \m{\EExt^1_{\ko_n}(\ke_1\ot\L^\vee,
\ko_n)=0} , on en d\'eduit la suite exacte
\[0\lra(\ke/\ke^{(1)})^{\vee_n}\lra\ke^{\vee_n}\lra\ke^{(1)*}\ot L^{n-1}\lra 0
.\]
De m\^eme en dualisant encore une fois on obtient la suite exacte
\[0\lra\ke^{(1)}\lra\ke^{\vee\vee}\lra\ke_1\ot\L^\vee\lra 0 \]
et le diagramme commutatif
\xmat{
0\ar[r] & \ke^{(1)}\ar[r]\fleq[d] & \ke\ar[r]\ar[d] &
\ke_1\ot\L^\vee\ar[r]\fleq[d] & 0\\
0\ar[r] & \ke^{(1)}\ar[r] & \ke^{\vee\vee}\ar[r] & \ke_1\ot\L^\vee\ar[r] & 0
}
Il en d\'ecoule que le morphisme canonique \ \m{\ke\to\ke^{\vee\vee}} \ est un
isomorphisme, c'est-\`a-dire que $\ke$ est r\'eflexif. Le fait que \
\m{\EExt^i_{\ko_n}(\ke,\ko_n)=0} \ pour \m{i\geq 1} d\'ecoule aussi de la
suite exacte (\ref{eq1}) et du fait que le th\'eor\`eme est vrai pour \m{n-1}.
On a donc montr\'e que (ii) entraine (i) et (iii).

Il reste \`a montrer que si $\ke$ est un faisceau coh\'erent sur \m{C_n} ayant
de la torsion, alors on a \ \m{\EExt^1_{\ko_n}(\ke,\ko_n)\not=0}. Soit $\T$
le sous-faisceau de torsion de $\ke$. Alors \m{\ke/\T} est sans torsion. En
utilisant la suite exacte
\[0\lra\T\lra\ke\lra\ke/\T\lra 0\]
et le fait que (ii) entraine (iii) on voit que
\[\EExt^1_{\ko_n}(\ke,\ko_n) \ \simeq \ \EExt^1_{\ko_n}(\T,\ko_n) .\]
Mais \m{\EExt^1_{\ko_n}(\T,\ko_n)} contient le faisceau non nul \
\m{\EExt^1_{\ko_n}(\T_{\mid C},\ko_n)=\widetilde{\T_{\mid C}}\ot L^{n-1}} ,
donc \Nligne \m{\EExt^1_{\ko_n}(\ke,\ko_n)\not=0}.
\end{proof}

\sepprop

\begin{subsub}{\bf Corollaire : }\label{cor1_refl} Soit \
\m{0\to\ke\to\kf\to\kg\to 0} \ une suite exacte de faisceaux coh\'erents sans
torsion sur \m{C_n}. Alors la suite duale \
\m{0\to\kg^\vee\to\kf^\vee\to\ke^\vee\to 0} \ est aussi exacte.
\end{subsub}

\sepprop

\begin{subsub}{\bf Corollaire : }\label{cor2_refl} Soit $\ke$ un faisceau
coh\'erent sur \m{C_n}. Alors on a \ \m{\EExt^i(\ke,\ko_n)=0} \ pour \m{i\geq
2}.
\end{subsub}

\begin{proof}
Il existe un morphisme surjectif \ \m{\E\to\ke} , $\E$ \'etant un fibr\'e
vectoriel sur \m{C_n} (on peut prendre pour $\E$ une somme directe de copies de
\m{\ko(-k)}, \m{\ko(1)} \'etant un fibr\'e en droites tr\`es ample sur \m{C_n}
et $k$ un entier assez grand). On a donc une suite exacte
\[0\lra\kf\lra\E\lra\ke\lra 0 ,\]
o\`u $\kf$ est un faisceau sans torsion, donc r\'eflexif. Le r\'esultat
d\'ecoule donc de la suite exacte longue obtenue en dualisant la pr\'ec\'edente,
et du th\'eor\`eme \ref{theo_refl}.
\end{proof}

\end{sub}

\sepsub

\Ssect{Dualit\'e de Serre pour les faisceaux r\'eflexifs}{Ser}

\begin{subsub}\label{du_se}{\bf Th\'eor\`eme : } Soit $\kf$ un faisceau
coh\'erent r\'eflexif sur \m{C_n}. Alors on a des isomorphismes fonctoriels
\[H^i(C_n,\kf) \ \simeq \ H^{1-i}(C_n,\kf^\vee\ot\omega_{C_n})^*\]
pour \m{i=0,1}.
\end{subsub}

\begin{proof}
Le th\'eor\`eme de dualit\'e de Serre (cf. par exemple \cite{ha}, theorem III,
7.6) donne les isomorphismes
\[H^i(C_n,\kf)\simeq\Ext^{1-i}_{\ko_n}(\kf,\omega_{C_n})^* .\]
D'apr\`es la suite spectrale des Ext et le th\'eor\`eme \ref{theo_refl} on a
\[\Ext^{1-i}_{\ko_n}(\kf,\omega_{C_n})\simeq
H^{1-i}(C_n,\kf^\vee\ot\omega_{C_n}) .\]
\end{proof}

\sepprop

\begin{subsub}\label{rem_se}{\bf Remarque : }\rm
La g\'en\'eralisation du r\'esultat pr\'ec\'edent aux Ext, c'est-\`a-dire le  
fait que
\[\Ext_{\ko_n}^i(\kf,\kg)\simeq\Ext_{\ko_n}^{1-i}(\kg,\kf\ot\omega_{C_n})^*\]
pour des faisceaux sans torsion $\kf$, $\kg$, est fausse en g\'en\'eral,
contrairement \`a ce qui se passe sur les vari\'et\'es lisses (cf. \cite{dr_lp},
prop. (1.2)). Par exemple, si \m{1<m<n}, la r\'esolution canonique de \m{\ko_i}
\[\cdots\lra\L^{2n}\lra\L^{n+m}\lra\L^n\lra\L^m\lra\ko_n\lra\ko_m\]
montre que
\[\Ext_{\ko_n}^{2i}(\ko_m,\ko_C)=L^{in} , \ \ \ \
\Ext_{\ko_n}^{2i+1}(\ko_m,\ko_C)=L^{in+m}\]
pour tout \m{i\geq 0}, ce qui exclue la dualit\'e pr\'ec\'edente, qui
entrainerait l'annulation des \m{\Ext^i} pour \m{i\geq 2}.
\end{subsub}

\end{sub}

\sepsec

\section{R\'esolutions p\'eriodiques des faisceaux quasi localement
libres}\label{res_qll}

On consid\`ere dans ce chapitre une courbe multiple primitive \m{C_n} de courbe
r\'eduite associ\'ee $C$. On utilise les notations de \ref{def_nota}.

\sepsub

\Ssect{Fibr\'es vectoriels et faisceaux quasi localement libres}{FV_qll}

\begin{subsub}\label{perio_pro1}{\bf Proposition : } Soit $\ke$ un faisceau
quasi localement libre sur \m{C_n}. Alors il existe un fibr\'e vectoriel $\E$
sur \m{C_n} et un morphisme surjectif \ \m{\E\to\ke} \ induisant un isomorphisme
\ \m{\E_{\mid C}\simeq\ke_{\mid C}} .
\end{subsub}

\begin{proof}
Par r\'ecurrence sur $n$. Le cas $n=1$ est \'evident. Supposons donc que \m{n>1}
et que le r\'esultat est vrai pour \m{n-1}. Si le support de $\ke$ est contenu
dans \m{C_{n-1}} il existe un fibr\'e vectoriel $\F$ sur \m{C_{n-1}} et un
morphisme surjectif \m{\F\to\ke} induisant un isomorphisme \ \m{\F_{\mid
C}\simeq\ke_{\mid C}}. D'apr\'es \cite{dr2}, th\'eor\`eme 3.1.1, on peut
\'etendre $\F$ en un fibr\'e vectoriel $\E$ sur \m{C_n}, et on obtient ainsi le
r\'esultat voulu. On peut donc supposer que le support de $\ke$ est \m{C_n}.
Le morphisme quotient \ \m{\ke\to\ke/\ke_{n-1}} induit un isomorphisme \
\m{\ke_{\mid C}\simeq(\ke/\ke_{n-1})_{\mid C}}, et d'apr\`es l'hypoth\`ese
de r\'ecurrence, puisque \m{\ke/\ke_{n-1}} est de support \m{C_{n-1}}, il existe
un fibr\'e vectoriel \m{\E_0} sur \m{C_{n-1}} et un morphisme surjectif \
\m{\phi_0:\E_0\to\ke/\ke_{n-1}} \ induisant un isomorphisme \ \m{\E_{0\mid
C}\simeq(\ke/\ke_{n-1})_{\mid C}}. On cherche un fibr\'e vectoriel $\E$ sur
\m{C_n} prolongement de \m{\E_0} tel qu'on ait un diagramme commutatif avec
lignes exactes
\xmat{0\ar[r] & \E_{\mid C}\ot L^{n-1}=\ke_{\mid C}\ot L^{n-1}\ar[r]\flon[d] &
\E\ar[r]\ar[d]^\phi & \E_0\ar[r]\flon[d]^{\phi_0} & 0\\
0\ar[r] & \ke_{n-1}\ar[r] & \ke\ar[r] & \ke/\ke_{n-1}\ar[r] & 0
}
o\`u la fl\`eche verticale de gauche est le morphisme canonique. Le morphisme
$\phi$ est alors surjectif et induit un isomorphisme \m{\E_{\mid
C}\simeq\ke_{\mid C}}. 

Soit \m{\sigma\in\Ext^1_{\ko_n}(\ke/\ke_{n-1},\ke)}, correspondant \`a la suite
exacte du bas. Soit \m{\sigma'\in\Ext^1_{\ko_n}(\E_0,\ke_{n-1})} l'image de
$\sigma$ dans \m{\Ext^1_{\ko_n}(\E_0,\ke_{n-1})}. Il suffit de montrer qu'il
existe \m{\sigma''\in\Ext^1_{\ko_n}(\E_0,\ke_{\mid C}\ot L^{n-1})} tel que dans
l'extension associ\'ee
\[0\lra\ke_{\mid C}\ot L^{n-1}\lra\E\lra\E_0\lra 0 ,\]
$\E$ soit localement libre, et que l'image de \m{\sigma''} dans
\m{\Ext^1_{\ko_n}(\E_0,\ke_{n-1})} soit \m{\sigma'} (cf. \cite{dr3}, corollaire
4.3.3). D'apr\`es \ref{Ext_ii} on a un diagramme commutatif avec lignes exactes
\xmat{
\Ext^1_{\ko_C}(\ke_{\mid C},\ke_{\mid C}\ot L^{n-1})\flinc[r]^-{i_0}
\ar[d]^\gamma & \Ext^1_{\ko_n}(\E_0,\ke_{\mid C}\ot L^{n-1})\ar[r]\flon[d]^\rho
 & \End(\ke_{\mid C})\ar[d]^\theta \\
\Ext^1_{\ko_C}(\ke_{\mid C},\ke_{n-1})\flinc[r]^-i &
\Ext^1_{\ko_n}(\E_0,\ke_{n-1})\flon[r]^-\lambda & \Hom(\ke_{\mid C}\ot L^{n-1},
\ke_{n-1})
}
o\`u les fl\`eches verticales sont induites par le morphisme canonique \
\m{\mu:\ke_{\mid C}\ot L^{n-1}\to\ke_{n-1}} . On a \ \m{\lambda(\sigma')=\mu=
\theta(I_{\ke_{\mid C}})} . Soit \ \m{\sigma_0\in\Ext^1_{\ko_n}(\E_0,
\ke_{\mid C}\ot L^{n-1})} \ dont l'image dans \m{\End(\ke_{\mid C})} est
\m{I_{\ke_{\mid C}}}. Alors on a
\[\rho(\sigma_0)-\sigma' \ \in \ i(\Ext^1_{\ko_n}(\ke_{\mid C},\ke_{n-1})) .\]
Puisque $\mu$ est surjectif il en est de m\`eme de $\gamma$. Donc il existe
\m{\sigma_1\in\Ext^1_{\ko_n}(\ke_{\mid C},\ke_{\mid C}\ot L^{n-1})} tel que
\[\rho(\sigma_0)-\sigma' \ = \ i\circ\gamma(\sigma_1) \ = \ \rho\circ
i_0(\sigma_1) .\]
On a alors \ \m{\sigma''=\sigma_0-i_0(\sigma_1)}.
\end{proof}

\sepprop

Soit $\ke$ un faisceau quasi localement libre sur \m{C_n}. Soit \ \m{\phi:\E\to
\ke} \ un morphisme surjectif, $\E$ \'etant localement libre, induisant un
isomorphisme \m{\E_{\mid C}\simeq\ke_{\mid C}}. Le faisceau \m{\kn=\ker(\phi)}
est aussi quasi localement libre, d'apr\`es \cite{dr2}, th\'eor\`eme 5.2.1.

\sepprop

\begin{subsub}\label{perio_lem1}{\bf Lemme : } Soit \m{(m_1,\ldots,m_n)} le type
de $\ke$. Alors le type de $\kn$ est\Nligne \m{(m_{n-1},m_{n-2},\ldots,m_1,0)}.
\end{subsub}

\begin{proof}
Soit $P$ un point ferm\'e de $C$. On a \
\m{\ke_P\simeq\bigoplus_{i=1}^nm_i\ko_{iP}} . On en d\'eduit la suite \m{(e_j)}
constitu\'ee des \'el\'ements 1 des facteurs \m{\ko_{iP}}. Soit \m{\ov{e_j}}
l'image de \m{e_j} dans \m{\ke_{\mid C,P}}. Alors \m{(\ov{e_j})} est une base de
\m{\ke_{\mid C,P}}. Soit \m{f_j\in\E_P} au dessus de \m{e_j}, \m{\ov{f_j}}
l'image de \m{f_j} dans \m{\E_{\mid C,P}}. Alors \m{(\ov{f_j})} est une base de
\m{\E_{\mid C,P}}, puisque c'est l'image de la base \m{(\ov{e_j})} par
l'isomorphisme \m{\ke_{\mid C,P}\simeq\E_{\mid C,P}}. En utilisant
l'isomorphisme \ \m{\E_P\simeq\rg(\E)\ko_{nP}} \ d\'efini par cette base on
voit que \m{\phi_P} est la somme des morphismes canoniques
\m{\ko_{nP}\to\ko_{iP}}.
Le r\'esultat en d\'ecoule imm\'ediatement.
\end{proof}

\sepprop

\begin{subsub}\label{perio_rem}{\bf Remarque : } \rm Plus g\'en\'eralement, on
d\'emontre de la m\^eme fa\c con que pour tout morphisme surjectif
\m{f:\F\to\ke}, $\F$ \'etant un fibr\'e vectoriel de rang $r$ sur \m{C_n},
\m{\ker(f)} est un faisceau quasi localement libre de type
\m{(m_{n-1},m_{n-2},\ldots,m_1,r-\sigg_{i=1}^nm_i)}.
\end{subsub}

\sepprop

Il existe un lien \'etroit entre les faisceaux $\ke$ et $\kn$. Dans le
r\'esultat suivant on compare les  morphismes \m{\lambda_{ij}} et \m{\mu_{ij}}
pour $\ke$ et $\kn$ (cf. \ref{rel_filt}) :

\sepprop

\begin{subsub}\label{perio_pro2}{\bf Proposition : } Soient $q$, $k$ des entiers
tels que \m{0<k<q<n}. Alors il existe des isomorphismes canoniques
\[\coker(\lambda_{q-1,k}(\kn))\simeq\ker(\mu_{n-q,k}(\ke)) , \quad\quad
\coker(\lambda_{q-1,k}(\ke))\ot L^n\simeq\ker(\mu_{n-q,k}(\kn)) .\]
\end{subsub}

\begin{proof} Soient $P$ un point ferm\'e de $C$ et \m{z\in\ko_{nP}} une
\'equation de $C$.

On d\'efinit d'abord le premier isomorphisme. Soient \
\m{u\in\kn_P^{(q-k)}/\kn_P^{(q-k-1)}} \ et \ \m{\ov{u}\in\kn_P^{(q-k)}} \ au
dessus de $u$. On a \ \m{z^{q-k}\ov{u}=0}, donc \ \m{\ov{u}\in z^{n-q+k}\E_P} ,
et on peut \'ecrire $\ov{u}$ sous la forme \m{\ov{u}=z^{n-q+k}v}. On a
\m{\ov{u}\in\kn_P}, donc \ \m{z^{n-q+k}\phi_P(v)=0}. Soit \m{e\in\ke_{\mid C,P}}
l'image de \m{\phi_P(v)}. On a donc \
\m{\mu_{n-q,k,P}(z^{n-q}e\ot z^k)=0} . On voit ais\'ement que \ \m{
z^{n-q}e\ot z^k\in(\ke_{n-q}\ot L^k)_P} \ est ind\'ependant des choix
effectu\'es. On d\'efinit donc ainsi un morphisme
\[\ov{\psi}:\kn^{(q-k)}/\kn^{(q-k-1)}\lra\ker(\mu_{n-q,k}(\ke)) .\]
Il est ais\'e de voir que si \ \m{u\in\imm(\lambda_{q-1,k})_P} , alors
\m{\ov{\psi}_P(u)=0}. On obtient donc un morphisme
\[\psi:\coker(\lambda_{q-1,k}(\kn))\lra\ker(\mu_{n-q,k}(\ke)) .\]
On d\'efinit maintenant le morphisme inverse. Soit \m{w\in\ke_{n-q,k,P}} tel que
\m{\mu_{n-q,k,P}(w\ot z^k)=0} . Soit \m{\beta\in\ke_P} tel que \m{z^{n-q}\beta}
soit au dessus de $w$. Il existe alors \m{\gamma\in\ke_P} tel que \Nligne
\m{z^{n-q+k}\beta=z^{n-q+k+1}\gamma}, c'est-\`a-dire \ \m{z^{n-q+k}(\beta-
z\gamma)=0}. Puisque \m{z^{n-q}(\beta-z\gamma)=0} est aussi au dessus de $w$, on
peut supposer que \m{\gamma=0}, c'est-\`a-dire \m{z^{n-q+k}\beta=0}. Soit
\m{u\in\E_P} au dessus de $\beta$. On a donc \ \m{z^{n-q+k}u\in\kn^{(q-k)}_P} .
On d\'efinit le morphisme inverse en associant \`a \m{w\ot z^k} l'image de
\m{z^{n-q+k}u} dans \m{\coker(\mu_{n-q,k}(\kn)_P)}. Les v\'erifications sont
laiss\'ees au lecteur.

La d\'efinition du second isomorphisme est analogue.
\end{proof}

\sepprop

\begin{subsub}\label{perio_cor}{\bf Corollaire : } On a, pour \m{0\leq i < n},
\[\rg(G_i(\kn)) \ = \ \rg(G^{(1)}(\ke))-\rg(G^{(n-i)}(\ke)) ,\]
\[\deg(G_i(\kn)) \ = \ \deg(G^{(1)}(\ke))-\deg(G^{(n-i)}(\ke))
+\big((i+1)\rg(G^{(1)}(\ke))-n\rg(G^{(n-i)}(\ke))\big)\deg(L) .\]
\end{subsub}

\begin{proof} Posons
\[\rho_i=\rg(G_i(\kn)), \quad \delta_i=\deg(G_i(\kn)), \quad
s_i=\rg(G^{(i+1)}(\ke)), \quad e_i=\deg(G^{(i+1)}(\ke)) .\]
On d\'eduit du
second isomorphisme de la proposition \ref{perio_pro2} les \'egalit\'es
\[\rho_{n-q}-\rho_{n-q+1} \ = \ s_{q-2}-s_{q-1} ,\]
\[\delta_{n-q}-\delta_{n-q+1} \ = \ e_{q-2}-e_{q-1}+\big(s_{q-2}n-s_{q-1}(n+1)
-\rho_{n-q}\big)\deg(L) .\]
Puisque $\kn$ est concentr\'e sur \m{C_{n-1}} on a \
\m{\rho_{n-1}=\delta_{n-1}=0} , et les \'egalit\'es pr\'ec\'edentes donnent
pour \m{q=2} le corollaire \ref{perio_cor} pour \m{i=n-2}. On en d\'eduit les
autres cas par r\'ecurrence descendante sur $i$.
\end{proof}

\end{sub}

\sepsub

\Ssect{Support des {\rm Ext} de faisceaux coh\'erents}{sup_ext}

\begin{subsub}\label{perio_th1} {\bf Th\'eor\`eme : } Soient $\ke$, $\kf$ des
faisceaux coh\'erents sur \m{C_n}, avec $\ke$ sans torsion. Alors le support du
faisceau \m{\EExt^1_{\ko_n}(\ke,\kf)} est contenu dans \m{C_{n-1}}.
\end{subsub}

\begin{proof}
Soient $P$ un point ferm\'e de $C$, \m{z\in\ko_{nP}} une \'equation de $C$, et
\m{x\in\ko_{nP}} au dessus d'une \'equation de $P$ dans \m{\ko_{CP}}. Il suffit
de prouver l'assertion locale analogue suivante : soient $M$, $N$ des
\m{\ko_{nP}}-modules de type fini, $M$ \'etant sans torsion. Alors on a \
\m{z^{n-1}.\Ext^1_{\ko_{nP}}(M,N)=\nsp} . Pour tout \m{\ko_{nP}}-module $W$, on
note \m{X_W:W\to W} la multiplication par \m{z^{n-1}}.

Soit $\E$ un \m{\ko_{nP}}-module libre de type fini tel qu'on ait un morphisme
surjectif \ \m{\pi:\E\to M} . Soit \m{V=\ker(\pi)}. On a donc une suite exacte \
\m{0\to V\to\E\to M\to 0} , d'o\`u on d\'eduit la suivante :
\xmat{\Hom(\E,N)\ar[r]^-{r_N} & \Hom(V,N)\ar[r] & \Ext^1_{\ko_{nP}}(M,N)\ar[r]
& 0 .}
Il faut donc montrer que si \m{\alpha\in\Hom(V,N)}, alors \ \m{\alpha\circ
X_V\in\imm(r_N)}. On a un diagramme commutatif
\xmat{\Hom(\E,V)\ar[r]^-{r_V}\ar[d]^{\Hom(\E,\alpha)} &
\Hom(V,V)\ar[d]^{\Hom(V,\alpha)} \\ \Hom(\E,N)\ar[r]^-{r_N} & \Hom(V,N)}
Supposons que l'assertion soit vraie pour \m{N=V} et \m{\alpha=I_V},
c'est-\`a-dire qu'il existe\Nligne \m{\lambda\in\Hom(\E,V)} tel que
\m{r_V(\lambda)=X_V}. On a alors en g\'en\'eral
\[\alpha\circ X_V=\alpha\circ r_V(\lambda)=r_N(\alpha\circ\lambda) ,\]
d'apr\`es le diagramme commutatif pr\'ec\'edent. Il suffit donc de traiter le
cas o\`u \m{N=V} et \m{\alpha=I_V}, c'est-\`a-dire montrer que \m{X_V} se
prolonge en un morphisme \m{\E\to V}.

Posons \ \m{V\ot\ko_{CP}=V_0\oplus T_0}, o\`u \m{T_0} est de torsion et \m{V_0}
est un \m{\ko_{CP}}-module libre. Il existe des bases \m{(e_i)_{1\leq i\leq
\rg(V_0)}}, \m{(f_i)_{1\leq i\leq \rg(\E)}} de \m{V_0}, $\E$ respectivement, et
un entier $p$, \m{0\leq p\leq \rg(V_0)}, tels que le morphisme \
\m{\psi:V_0\to\E\ot\ko_{CP}} \ d\'eduit de l'inclusion \m{V\subset\E} soit tel
que \ \m{\psi(e_i)=x^{k_i}f_i} \ pour \m{1\leq i\leq p}, avec \m{k_i\geq 0}, et
\m{\psi(e_i)\in z\E} si \m{i>p}. Soient \m{\epsilon_i\in V}, \m{\phi_j\in\E} au
dessus de \m{e_i}, \m{f_j} respectivement. Alors on a (en voyant $V$ comme un
sous-module de $\E$), \m{\epsilon_i=x^{k_i}\phi_i+zu_i}, avec \m{u_i\in\E}. Donc
\ \m{z^{n-1}\epsilon_i=x^{k_i}z^{n-1}\phi_i} . Puisque $M$ est sans torsion, on
a \ \m{z^{n-1}\phi_i\in V} . On d\'efinit maintenant \ \m{\Theta:\E\to V} \ par
: \m{\Theta(\phi_i)=z^{n-1}\phi_i} \ si \m{1\leq i\leq p}, et
\m{\Theta(\phi_i)=0} si \m{i>p}. On v\'erifie ais\'ement que \ \m{\Theta_{\mid
V}=X_V} .
\end{proof}

\sepprop

\begin{subsub}\label{perio_rem1}{\bf Remarque : }\rm L'hypoth\`ese que $\ke$ est
sans torsion est n\'ecessaire dans le th\'eor\`eme pr\'ec\'edent. Par exemple si
$\T$ est un faisceau de torsion sur \m{C_n} non contenu dans \m{C_{n-1}}, le
support de \m{\widetilde{\T}=\EExt^1_{\ko_n}(\T,\ko_n)}  n'est pas contenu dans
\m{C_{n-1}}.
\end{subsub}

\sepprop

\begin{subsub}\label{perio_cor1} {\bf Corollaire : } Soient $\ke$, $\kf$ des
faisceaux coh\'erents sur \m{C_n}. Alors, pour tout entier \m{i\geq 2}, le
support du faisceau \m{\EExt^i_{\ko_n}(\ke,\kf)} est contenu dans \m{C_{n-1}}.
\end{subsub}

\begin{proof}
Soit $\E$ un fibr\'e vectoriel sur \m{C_n} tel qu'on ait un morphisme surjectif
\m{\phi:\E\to\ke}. Soit \m{\kn=\ker(\Phi)}. On a donc une suite exacte \
\m{0\to\kn\to\E\to\ke\to 0}. On en d\'eduit que \
\m{\EExt^2_{\ko_n}(\ke,\kf)\simeq\EExt^1_{\ko_n}(\kn,\kf)} . D'apr\`es le
th\'eor\`eme \ref{perio_th1}, le support de \m{\EExt^2_{\ko_n}(\ke,\kf)} est
donc contenu dans \m{C_{n-1}}. Le cas des \m{\EExt^i} pour \m{i\geq 2} se traite
par r\'ecurrence sur $i$.
\end{proof}

\sepprop

\begin{subsub}\label{perio_cor2} {\bf Corollaire : } Soient $\ke$, $\kf$ des
faisceaux coh\'erents sur \m{C_n}, \m{\D} un fibr\'e en droites sur \m{C_n} tel
que \ \m{\D_{\mid C_{n-1}}\simeq\ko_{n-1}} . Alors, pour tout entier \m{i\geq
2}, il existe un isomorphisme canonique et fonctoriel \
\m{\Ext^i_{\ko_n}(\ke,\kf\ot\D)\simeq\Ext^i_{\ko_n}(\ke,\kf)} .
\end{subsub}

\begin{proof} Pour tout entier \m{j\geq 2} on a
\[\EExt^j_{\ko_n}(\ke,\kf\ot\D) \ \simeq \ \EExt^j_{\ko_n}(\ke,\kf)\ot\D \
\simeq\EExt^j_{\ko_n}(\ke,\kf)\]
d'apr\`es le corollaire \ref{perio_cor1}, donc le r\'esultat d\'ecoule de la
suite spectrale des Ext et de sa fonctorialit\'e.
\end{proof}

\end{sub}

\sepsub

\Ssect{Morphismes et classes canoniques}{mo_ca}

Soit $N$ un entier tel que \m{N\gg n}. D'apr\`es \cite{dr1} il existe une courbe
multiple primitive \m{C_N} extension de \m{C_n} tel que le faisceau d'id\'eaux
de $C$ dans \m{C_N} restreint \`a \m{C_n} soit isomorphe \`a $\L$. On peut donc
noter aussi $\L$ un fibr\'e en droites sur \m{C_N} extension de \m{\ki_C}.

\sepprop

\begin{subsub}\label{perio_lem2}{\bf Lemme : } Soit $\ke$ un faisceau coh\'erent
sur \m{C_n}. Alors on a un isomorphisme canonique \
\m{\Tor^1_{\ko_N}(\ke,\ko_n)\simeq\ke\ot\L^n} .
\end{subsub}

\begin{proof}
Cela d\'ecoule de la r\'esolution localement libre de \m{\ko_n} sur \m{C_N}
\[\ldots\lra\L^N\lra\L^n\lra\ko_N\lra\ko_n\lra 0\]
et du fait que \m{N\gg n}.
\end{proof}

\sepprop

Soient $\ke$, $\kf$ des faisceaux coh\'erents sur \m{C_n}. D'apr\`es la
proposition \ref{extf_pr} on a une suite exacte
\xmat{0\ar[r] & \Ext^1_{\ko_n}(\ke,\kf\ot\L^n)\ar[r] &
\Ext^1_{\ko_N}(\ke,\kf\ot\L^n)\ar[r] & \Hom(\Tor^1_{\ko_N}(
\ke,\ko_n),\kf\ot\L^n)\\
\ar[r]^-{\Lambda_{\ke\kf}} & \Ext^2_{\ko_n}(\ke,\kf\ot\L^n)}
On a donc d'apr\`es le lemme \ref{perio_lem2} un morphisme canonique fonctoriel
\[\Lambda_{\ke\kf}:\Hom(\ke,\kf)\lra\Ext^2_{\ko_n}(\ke,\kf\ot\L^n) .\]
D'apr\`es le corollaire \ref{perio_cor2}, \m{\Ext^2_{\ko_n}(\ke,\kf\ot\L^n)} est
ind\'ependant du choix de $\L$.
En utilisant la preuve dans \cite{dr2} de la proposition \ref{extf_pr} on montre
ais\'ement que \m{\Lambda_{\ke\kf}} est ind\'ependant du plongement de \m{C_n}
dans \m{C_N} et est fonctoriel par rapport \`a $\ke$ et $\kf$. Les morphismes
\m{\Lambda_{\ke\kf}} sont donc enti\`erement d\'etermin\'es par les classes
\[\lambda_\ke \ = \ \Lambda_{\ke\ke}(I_\ke) \ \in \
\Ext^2_{\ko_n}(\ke,\ke\ot\L^n) .\]



\end{sub}

\sepsub

\Ssect{R\'esolutions p\'eriodiques}{re_pe}

On utilise comme dans \ref{mo_ca} un plongement \m{C_n\subset C_N} avec \m{N\gg
n}. Soit $\ke$ un faisceau quasi localement libre sur \m{C_n}. D'apr\`es la
proposition \ref{perio_pro1} il existe un fibr\'e vectoriel $\V$ sur \m{C_N} et
un morphisme surjectif \m{\phi:\V\to\ke} induisant un isomorphisme \m{\V_{\mid
C}\simeq\ke_{\mid C}}. Soit \m{\U=\ker(\phi)}. On a donc une suite exacte sur
\m{C_N} :
\[0\lra\U\lra\V\lra\ke\lra 0 .\]
Soient \m{\E=\V_{\mid C_n}}, \m{\F=\U_{\mid C_n}} . Il d\'ecoule du lemme
\ref{perio_lem1} et du fait que \m{N\gg n} que $\F$ est un fibr\'e vectoriel sur
\m{C_n}, de m\^eme rang que $\E$. En restreignant la suite exacte pr\'ec\'edente
\`a \m{C_n} et en utilisant le lemme \ref{perio_lem2} on obtient la suite exacte
sur \m{C_n}
\[0\lra\ke\ot\L^n\lra\F\lra\E\lra\ke\lra 0 .\]
En utilisant la d\'emonstration de la proposition \ref{extf_pr} on montre
ais\'ement que l'\'el\'ement de \m{\Ext^2_{\ko_n}(\ke,\ke\ot\L^n)} associ\'e \`a
la suite exacte pr\'ec\'edente est \m{\lambda_\ke}. On en d\'eduit le

\sepprop

\begin{subsub}\label{perio_th2}{\bf Th\'eor\`eme : } Soient $\ke$, $\kf$ des
faisceaux coh\'erents, avec $\ke$ quasi localement libre.

1 - La multiplication par \m{\lambda_\ke}
\[\HHom(\ke,\kf)\lra\EExt^2_{\ko_n}(\ke,\kf)\ot\L^n\]
est surjective.

2 - Pour tout entier \m{i\geq 1} la multiplication par \m{\lambda_\ke}
\[\EExt^i_{\ko_n}(\ke,\kf)\lra\EExt^{i+2}_{\ko_n}(\ke,\kf)\ot\L^n\]
est un isomorphisme.
\end{subsub}

\sepprop

\begin{subsub}\label{perio_cor4}{\bf Corollaire : } Soit $\ke$ un faisceau
quasi localement libre sur \m{C_n}. Alors $\ke$ est localement libre si et
seulement si on a \ \m{\lambda_\ke=0}.
\end{subsub}

\begin{proof}
Si $\ke$ est localement libre, on a
\[\Ext^2_{\ko_n}(\ke,\ke\ot\L^n)=H^2(\ke^\vee\ot\ke\ot\L^n)=\nsp \quad,\]
donc \ \m{\lambda_\ke=0}. R\'eciproquement, supposons que \ \m{\lambda_\ke=0}.
D'apr\`es le th\'eor\`eme \ref{perio_th2}, 2-, on a
\[\EExt^1_{\ko_n}(\ke,\ko_C) \ = \ 0 .\]
Supposons que $\ke$ ne soit pas localement libre. Alors $\ke$ est localement
isomorphe \`a une somme directe de faisceaux du type \m{\ko_k} dont l'un
d'entre eux au moins est tel que \m{k<n}. Il suffit donc de montrer que si $P$
est un point ferm\'e de $C$, alors on a
\[\Ext^1_{\ko_{nP}}(\ko_{kP},\ko_{CP}) \ \not= \ \nsp\quad .\]
Cela se voit ais\'ement en utilisant la r\'esolution canonique de \m{\ko_{kP}} :
\xmat{\cdots\ar[r] & \ko_{nP}\ar[r]^{\times z^k} & \ko_{nP}\ar[r]^{\times
z^{n-k}} & \ko_{nP}\ar[r]^{\times z^k} & \ko_{nP}\ar[r] & \ko_{kP}
}
(o\`u \ \m{z\in\ko_{nP}} \ est une \'equation de $C$).
\end{proof}

\sepprop

\begin{subsub}\label{perio_cor3}{\bf Corollaire : } Soient $\ke$, $\kf$ des
faisceaux coh\'erents, avec $\ke$ quasi localement libre.

1 - La multiplication par \m{\lambda_\ke}
\[\Ext^1_{\ko_n}(\ke,\kf)\lra\Ext^3_{\ko_n}(\ke,\kf\ot\L^n)\]
est surjective.

2 - Pour tout entier \m{i\geq 2} la multiplication par \m{\lambda_\ke}
\[\Ext^i_{\ko_n}(\ke,\kf)\lra\Ext^{i+2}_{\ko_n}(\ke,\kf\ot\L^n)\]
est un isomorphisme.
\end{subsub}

\end{sub}

\sepsec
\newpage
\section{Familles de faisceaux quasi localement libres}\label{fam_qll}

On consid\`ere dans ce chapitre une courbe multiple primitive \m{C_n} de courbe
r\'eduite associ\'ee $C$. On utilise les notations de \ref{def_nota}.

\sepsub

\Ssect{\'Etude locale des familles de faisceaux de type constant}{fam_qll_loc}

Soit $Y$ une vari\'et\'e alg\'ebrique int\`egre de dimension \m{d>0}. Soit $s$
un entier tel que \m{1\leq s\leq n}.
Soit $U$ un ouvert non vide de \m{Y\times C_n} tel que la projection de $U$
sur $Y$ soit surjective. Soit $\ke$ un faisceau coh\'erent sur $U$, plat sur
$Y$. On suppose que pour tout point ferm\'e \m{y\in Y}, \m{\ke_y} est de type
\m{(m_1^y,\ldots,m_n^y)}, avec \m{m_1^y\ldots,m_s^y} ind\'ependants de $y$, et
que pour \m{0\leq i<s}, \m{G_i(\ke_y)} est localement libre sur $C$.

On dit qu'un tel faisceau est une {\em famille de faisceaux de type constant
\`a l'ordre $s$} sur des ouverts de \m{C_n} (resp. sur
\m{C_n} si \ \m{U=Y\times C_n}) param\'etr\'ee par $Y$. Si \m{s=n} on dit que
$\ke$ est une {\em famille de faisceaux quasi localement libres de type
constant} sur des ouverts de \m{C_n} (resp. sur \m{C_n} si \ \m{U=Y\times C_n})
param\'etr\'ee par $Y$.

On d\'efinit comme pour les faisceaux coh\'erents sur \m{C_n} les filtrations
canoniques de $\ke$, c'est-\`a-dire les sous-faisceaux \m{\ke_i}, \m{\ke^{(i)}},
\m{1\leq i<n}.

Soient \m{P=(y,x)} un point ferm\'e de \m{Y\times C_n}, \m{z\in\ko_{nx}} une
\'equation de $C$ et $I$ l'id\'eal de \m{\lbrace y\rbrace\times C_n} dans
\m{{\bf m}=\ko_{ Y\times C_n,P}}.

\sepprop

\begin{subsub}\label{f_ql_lem1}{\bf Lemme : } Le faisceau \m{\ke_{\mid
U\cap(Y\times C)}} est localement libre de rang \m{m_1+\cdots+m_n}.
\end{subsub}

\begin{proof}
Cela d\'ecoule du fait que pour tout point ferm\'e $P$ de \m{Y\times C}, si
\m{m_P} d\'esigne l'id\'eal maximal de \m{\ko_{Y\times C_n,P}}, la dimension du
$\C$-espace vectoriel \ \m{\ke_P\ot(\ko_{Y\times C_n,P}/m_P)} \ est
\m{m_1+\cdots+m_n}, et du fait que \m{Y\times C} est r\'eduite.
\end{proof}

\sepprop

On pose \ \m{\V=\ke_{\mid U\cap(Y\times C)}}.

\sepprop

\begin{subsub}\label{f_ql_pro1}{\bf Proposition : } Soient
\m{\tau_1,\ldots,\tau_k\in I} dont les images dans \m{{\bf m}/{\bf m}^2} sont
lin\'eairement ind\'ependantes.

1 - Soit \m{u\in\ke_P} \ tel que \
\m{zu\in(\tau_1,\ldots,\tau_k)\ke_P}. Alors il existe \
\m{v\in(\tau_1,\ldots,\tau_k)\ke_P} \ tel que \m{zu=zv}.

2 - Soit \m{u\in(\tau_1,\ldots,\tau_k)\ke_P} tel que \m{zu=0}. Alors on peut
\'ecrire $u$ sous la forme \Nligne \m{u=\tau_1w_1+\cdots+\tau_kw_k} , avec 
\m{w_i} tels que \m{zw_i=0}.
\end{subsub}

\begin{proof}
Les deux assertions se d\'emontrent par r\'ecurrence sur $k$. Supposons d'abord
que \m{k=1} et soit \m{u\in\ke_P} tel que \m{zu\in(\tau_1)\ke_P} :
\m{zu=\tau_1v}. En des points \m{P'=(y',x')} voisins g\'en\'eriques de $P$,
\m{\tau_1} est inversible, donc $v$ est divisible par $z$ dans \m{\ke_P'}. Il en
d\'ecoule que l'image de $v$ dans \m{\V_{P'}} (la fibre en \m{P'} du fibr\'e
vectoriel $\V$) est nulle. Donc $v$ est nul comme section locale de $\V$. Donc
$v$ est multiple de $z$ : \m{v=zw}. On a donc \m{zu=z.\tau_1w}, ce qui
d\'emontre 1- pour \m{k=1}.

Soit \m{u\in(\tau_1)\ke_P} tel que \m{zu=0}. Posons \m{u=\tau_1v}. On a donc
\m{\tau_1.zv=0}. Puisque $\ke$ est plat sur $Y$, la multiplication par
\m{\tau_1} est injective. On a donc \m{zv=0}, ce qui d\'emontre 2- pour \m{k=1}.

Supposons que 1- soit vraie pour \m{k-1\geq 1}. Soit \m{u\in\ke_P}
tel que \m{zu\in(\tau_1,\ldots,\tau_k)\ke_P}. On consid\`ere la
sous-vari\'et\'e int\`egre d'un voisinage de $P$ dans $Y$, d\'efinie par
l'\'equation \m{\tau_k=0}. Alors \m{\ke_{\mid U\cap(Y'\times C_n)}} est plat sur
\m{Y'}. D'apr\`es l'hypoth\`ese de r\'ecurrence on peut \'ecrire \
\m{zu=zw+\tau_d\lambda}, avec \m{w\in(\tau_1,\ldots,\tau_{d-1})\ke_P},
c'est-\`a-dire \m{z(u-w)=\tau_d\lambda}. En faisant le m\^eme raisonnement que
dans le cas \m{k=1} on obtient que $\lambda$ est multiple de $z$ :
\m{\lambda=z\beta}. D'o\`u \ \m{zu=z(v+\tau_d\beta)}, et \
\m{v+\tau_d\beta\in(\tau_1,\ldots,\tau_k)\ke_P}, ce qui d\'emontre 1-.

Supposons que 2- soit vraie pour \m{k-1\geq 1}. Soit 
\m{u\in(\tau_1,\ldots,\tau_k)\ke_P} tel que \m{zu=0}. On consid\`ere la courbe
int\`egre \m{Y'\subset Y} d\'efinie au voisinage de $y$ par les \'equations
\m{\tau_2=\cdots=\tau_k=0}. Posons \ \m{u=\tau_1v_1+\cdots+\tau_kv_k}. Sur
\m{Y'} on a \ \m{u=\tau_1v_1}, et d'apr\`es le cas \m{k=1} on a \Nligne
\m{zv_1=0 {\rm\ mod }(\tau_2,\ldots,\tau_k)}. D'apr\`es 1- on peut \'ecrire
\m{zv_1=zw}, avec \m{w\in(\tau_2,\ldots,\tau_k)\ke_P} :\Nligne
\m{w=\tau_2w_2+\cdots+\tau_kw_k}. On a
\[z(u-\tau_1(v_1-w)) = \tau_2(v_2+t_1w_2)+\cdots+\tau_k(v_k+t_1w_k) ,\]
et \ \m{z(u-\tau_1(v_1-w))=0} . D'apr\`es l'hypoth\`ese de r\'ecurrence, on peut
\'ecrire
\[u-\tau_1(v_1-w)=\tau_2s_2+\cdots+\tau_ks_k ,\]
avec \m{zs_i=0} pour \m{2\leq i\leq k}, ce qui d\'emontre 2-.
\end{proof}

\sepprop

\begin{subsub}\label{f_ql_cor1}{\bf Corollaire : } Soit $i$ un entier tel que
\m{1\leq i<s}.

1 - Pour tout point ferm\'e $y$ de $Y$, les morphismes canoniques
\[(\ke_i)_y\to(\ke_y)_i , \quad\quad (\ke^{(i)})_y\to(\ke_y)^{(i)},
\quad\quad (\ke/\ke_i)_y\to\ke_y/(\ke_y)_i\]
sont des isomorphismes.

2 - Les faisceaux \m{\ke_i}, \m{\ke^{(i)}}, \m{\ke/\ke_i} sont plats sur $Y$,
et \m{\ke/\ke_i} est une famille de faisceaux quasi localement libres de type
constant \`a l'ordre $s$.
Les faisceaux \m{G_i(\ke)}, \m{G^{(i)}(\ke)} sont des fibr\'es vectoriels sur
\m{Y\times C}.

3 - Les degr\'es des faisceaux \m{(\ke_y)_i}, \m{(\ke_y)^{(i)}},
\m{\ke_y/(\ke_y)_i}, \m{G_i(\ke_y)}, \m{G^{(i)}(\ke_y)} sont ind\'ependants du
point ferm\'e \m{y\in Y}. En particulier, si \m{s=n}, les faisceaux \m{\ke_y}
sont de type complet constant (cf. \ref{type_complet}).
\end{subsub}

\begin{proof}
On montre d'abord que \ \m{(\ke^{(1)})_y\simeq(\ke_y)^{(1)}} . Soit \
\m{\theta:(\ke^{(1)})_y\to(\ke_y)^{(1)}} \ le morphisme canonique. Soient
\m{x\in C} un point ferm\'e et \m{P=(y,x)}. Soient \m{u\in\lbrack(\ke^{(1)})_y
\rbrack_x} et \m{\ov{u}\in\ke^{(1)}_P} au dessus de $u$. L'image $v$ de $u$
dans \m{(\ke_y)_x} est annul\'ee par $z$, donc \m{v\in\lbrack(\ke_y)^{(1)}
\rbrack_x}, et \m{\theta_x(u)=v}.

Montrons que \m{\theta_x} est injectif. Supposons que \m{\theta_x(u)=0}. On a
\m{z\ov{u}=0} et \m{\ov{u}\in I\ke_P}. D'apr\`es la proposition \ref{f_ql_pro1},
2-, on a \m{\ov{u}\in I(\ke^{(1)})_P}, donc \m{u=0} et \m{\theta_x} est
injectif.

Montrons que \m{\theta_x} est surjectif. Soit \m{w\in\lbrack(\ke_y)^{(1)}
\rbrack_x}. On a donc \m{zw=0}. Soit \m{\ov{w}\in\ke_P} au dessus de $w$. On a
\m{z\ov{w}\in I\ke_P}. D'apr\`es la proposition \ref{f_ql_pro1}, 1-, il existe
\m{\ov{v}\in I\ke_P} tel que \m{z\ov{w}=z\ov{v}}. On a donc \m{\ov{w}-\ov{v}
\in(\ke^{(1)})_P}, et si \m{\beta\in\lbrack(\ke^{(1)})_y\rbrack_x} est l'image
de \m{\ov{w}-\ov{v}}, on a \m{\theta_x(\beta)=w}. Donc \m{\theta_x} est
surjectif.

L'isomorphisme \m{(\ke/\ke_1)_y\simeq\ke_y/(\ke_y)_1} se d\'emontre de la
m\^eme fa\c con. L'isomorphisme\Nligne \m{(\ke_1)_y\simeq(\ke_y)_1} se d\'eduit
ais\'ement des deux pr\'ec\'edents. 

On d\'eduit de ce qui pr\'ec\`ede et de \cite{sga1}, expos\'e IV, cor. 5.7, que
\m{\ke_1} est plat sur $Y$.

Les cas \m{i>1} se d\'emontrent \`a partir du cas \m{i=1} par r\'ecurrence sur
$i$.
\end{proof}

\sepprop

Le r\'esultat suivant est une version relative de la proposition
\ref{perio_pro1} :

\sepprop

\begin{subsub}\label{f_ql_cor2}{\bf Corollaire : } Soient $X$ une vari\'et\'e
alg\'ebrique affine int\`egre et $\ke$ une famille de faisceaux quasi localement
libres de type constant sur \m{C_n} param\'etr\'ee par $X$. Alors il existe un
fibr\'e vectoriel $\E$ sur \m{X\times C_n} et un morphisme surjectif \
\m{\E\to\ke} \ induisant un isomorphisme
\ \m{\E_{\mid X\times C}\simeq\ke_{\mid X\times C}} .
\end{subsub}

\begin{proof} Analogue \`a la d\'emonstration de la proposition
\ref{perio_pro1}. On utilise le fait que \m{\ke_{\mid C}} est un fibr\'e
vectoriel sur \m{X\times C} (d'apr\`es le corollaire \ref{f_ql_cor1}), et que
pour tout faisceau coh\'erent $\kf$ sur \m{X\times C_n} on a \
\m{H^2(\kf)=\nsp} (cf. le d\'ebut de la d\'emonstration du th\'eor\`eme
\ref{prol_th}).
\end{proof}

\sepprop

\begin{subsub}\label{f_ql_th1}{\bf Th\'eor\`eme : } Soit $\ke$ une famille de
faisceaux quasi localement libres de type constant sur des ouverts de \m{C_n}
param\'etr\'ee par une vari\'et\'e int\`egre $Y$. Alors pour tout point ferm\'e
$P$ de \m{Y\times C_n} il existe un voisinage $V$ de $P$ dans $U$ tel que
\[\ke_{\mid V} \ \simeq \ \bigoplus_{i=1}^nm_i.p_{C_n}^*(\ko_i)_{\mid V} .\]
\end{subsub}

\begin{proof} On proc\`ede par r\'ecurrence sur $n$, le r\'esultat \'etant
\'evident si \m{n=1} (car $Y$ est int\`egre). Supposons que le r\'esultat est
vrai pour \m{n-1\geq 1}. On fait maintenant une r\'ecurrence sur \m{m_n}. Si
\m{m_n=0} le r\'esultat est vrai car alors $\ke$ est concentr\'e sur
\m{U\cap(Y\times C_{n-1})}. Supposons que \m{m_n>0} et que le r\'esultat est
vrai pour \m{m_n-1}.

Posons \m{P=(y,x)}. Pour tout ouvert $V$ de \m{Y\times C_n} et tout \m{y'\in
Y}, on note \m{V_{y'}=V\cap(\lbrace y'\rbrace\times C_n)}. On peut supposer que
$U$ est affine, que \ \m{\ke_{y\mid U_y}\simeq\oplus_{i=1}^nm_i\ko_{i\mid U_y}}
\ et que \ \m{\kl=p_{C_n}^*(\ki_C)} \ est un fibr\'e en droites trivial sur \
\m{Y\times C_{n-1}}. Soit \m{z\in H^0(\kl)} une section engendrant $\kl$. Soient
\ \m{\sigma\in H^0(\ke_y)} \ d\'efini par un \'el\'ement non nul de
\m{\C^{m_n}} et \m{\ov{\sigma}\in H^0(\ke)} un prolongement de $\sigma$. Alors \
\m{s=z^{n-1}\ov{\sigma}\in H^0(\ke_{n-1})} . D'apr\`es le corollaire
\ref{f_ql_cor1} \m{\ke_{n-1}} est un fibr\'e vectoriel sur \m{Y\times C},
et \ \m{s_{\mid U_y}\in H^0((\ke_y)_{n-1})} \ ne s'annule en aucun point de
\m{U_y}. Soit \m{V\in Y\times C} l'ouvert des points o\`u $s$ ne s'annule pas.
Soient \m{y'\in Y}, et \m{W\in V_{y'}} un ouvert tel que \ \m{\ke_{\mid W}
\simeq\oplus_{i=1}^nm_i\ko_{i\mid W}} . Le morphisme \ \m{\sigma_{\mid W} :
\ko_{n\mid W}\lra\ke_W} \ induit une section de \m{m_n\ko_C} qui ne s'annule en
aucun point. Il en d\'ecoule que
\[\coker(\sigma_{\mid W}) \ \simeq \biggl(\bigoplus_{i=1}^{n-1}m_i\ko_{i\mid
W}\biggr)\oplus (m_n-1)\ko_{n\mid W} .\]
Il d\'ecoule du corollaire 5.7 de \cite{sga1}, expos\'e IV que \
\m{\kf=\coker(\ov{\sigma}_{\mid V})} \ est plat sur $Y$. C'est une famille de
faisceaux quasi localement libres de type \m{(m_1,\ldots,m_{n-1},m_n-1)}.
D'apr\`es l'hypoth\`ese de r\'ecurrence on peut, quitte \`a remplacer $V$ par
un voisinage plus petit de $P$, supposer que
\[\kf \ \simeq \ \biggl(\bigoplus_{i=1}^{n-1}m_ip_{C_n}^*(\ko_{i})_{\mid
V}\biggr)\oplus (m_n-1)p_{C_n}^*(\ko_{n\mid W})_{\mid V} .\]
On a donc une suite exacte
\[0\lra\ko_V\lra\ke_V\lra\biggl(\bigoplus_{i=1}^{n-1}m_ip_{C_n}^*(\ko_{i})_{\mid
V}\biggr)\oplus (m_n-1)p_{C_n}^*(\ko_{n\mid W})_{\mid V} \lra 0 .\]
Mais on montre ais\'ement, en utilisant les r\'esolutions localement
libres habituelles des \m{\ko_i} sur \m{C_n}, que \ \m{\Ext^1_{\ko_V}(\kf,
\ko_V)=\nsp} . On obtient donc finalement \ \m{\ke_{\mid V}\simeq
\bigoplus_{i=1}^nm_i.p_{C_n}^*(\ko_i)_{\mid V}} .
\end{proof}

\sepprop

\begin{subsub}{\bf Remarque : }\label{f_ql_rem}\rm Le th\'eor\`eme
\ref{f_ql_th1} n'est pas vrai si on ne suppose pas que $Y$ est r\'eduite.
Supposons par exemple que $Y$ soit une courbe multiple primitive de
multiplicit\'e \m{m>0}, de courbe r\'eduite associ\'ee \m{C'}. On consid\`ere le
faisceau \ \m{\ke=p_Y^*(\ko_C')} \ sur \m{Y\times C_n}. C'est bien une famille
plate de faisceau coh\'erents sur \m{C_n} param\'etr\'ee par $Y$, on a \
\m{\ke_P\simeq\ko_n} \ pour tout point ferm\'e $P$ de $Y$. Mais la conclusion du
th\'eor\`eme \ref{f_ql_th1} est fausse pour cette famille de faisceaux.
\end{subsub}

\sepprop

\begin{subsub}{\bf Corollaire : }\label{f_ql_cor3} Soient $X$ une vari\'et\'e
alg\'ebrique affine int\`egre et $\ke$ une famille de faisceaux quasi localement
libres de type constant sur \m{C_n} param\'etr\'ee par $X$. Soient $x$ un point
ferm\'e de $X$. Alors le morphisme de d\'eformation infinit\'esimale de Koda\i
ira-Spencer de $\ke$ en $x$
\[\omega_x:TX_x\lra\Ext^1_{\ko_n}(\ke_x,\ke_x)\]
est \`a valeurs dans le sous-espace \m{H^1(\EEnd(\ke_x))} de
\m{\Ext^1_{\ko_n}(\ke_x,\ke_x)}.
\end{subsub}

\begin{proof} Pour tout point $q$ de $C$, la d\'eformation de \m{\ke_{(x,q)}}
induite par $\ke$ est triviale d'apr\`es la th\'eor\`eme \ref{f_ql_th1}. Il en
d\'ecoule que l'application compos\'ee
\xmat{TX_x\ar[r]^-{\omega_x} & \Ext^1_{\ko_n}(\ke_x,\ke_x)\ar[r] &
H^0(\EExt^1_{\ko_n}(\ke_x,\ke_x))}
est nulle. Donc \m{\omega_x} est \`a valeurs dans le noyau de l'application
canonique
\[\Ext^1_{\ko_n}(\ke_x,\ke_x)\to H^0(\EExt^1_{\ko_n}(\ke_x,\ke_x)) ,\]
qui est \m{H^1(\EEnd(\ke_x))}.
\end{proof}

\end{sub}

\sepsub

\Ssect{Irr\'eductibilit\'e}{irred}

Soit $\kp$ un ensemble de classes d'isomorphisme de faisceaux coh\'erents sur
\m{C_n}. On dit que $\kp$ est {\em irr\'eductible} si pour tous
\m{E_0,E_1\in\kp} il existe une famille plate $\ke$ de faisceaux coh\'erents
sur \m{C_n} param\'etr\'ee par une vari\'et\'e alg\'ebrique lisse
irr\'eductible $S$ telle que pour tout point ferm\'e $s$ de $S$ on ait
\m{\ke_s\in\kp}, et qu'il existe des points ferm\'es \m{s_0,s_1\in S} tels que \
\m{\ke_{s_0}\simeq E_0} et \m{\ke_{s_1}\simeq E_1}. 

On peut faire une d\'efinition semblable concernant des morphismes de faisceaux
coh\'erents sur \m{C_n}. On dit que deux morphismes \ \m{f:E\to F}, \m{f':E'\to
F'} \ sont {\it isomorphes} s'il existe des isomorphismes \ \m{\epsilon:E\to
E'}, \m{\phi:F\to F'} \ tels que le carr\'e
\xmat{E\ar[r]^f\ar[d]^\epsilon & F\ar[d]^\phi \\ E'\ar[r]^-{f'} & F'}
soit commutatif. Soit $\ks$ un ensemble de classes d'isomorphisme de morphismes
de faisceaux coh\'erents sur \m{C_n}. On dit que $\ks$ est {\em irr\'eductible}
si pour tous \m{f_0}, \m{f_1} dans $\ks$ il existe des familles plates $\ke$,
$\kf$ de faisceaux coh\'erents sur \m{C_n} param\'etr\'ees par une vari\'et\'e
alg\'ebrique lisse irr\'eductible $S$, et un morphisme \m{\theta:\ke\to\kf},
tels que pour tout point ferm\'e $s$ de $S$ le morphisme \m{\theta_s} soit dans
$\ks$, et qu'il existe des points ferm\'es \m{s_0,s_1} de $S$ tels que
\m{\theta_{s_0}=f_0} et \m{\theta_{s_1}=f_1}.

\sepprop

\begin{subsub}\label{irred_theo}{\bf Th\'eor\`eme : } Soient $r$, $d$,
\m{r_0,\ldots,r_{n-1}}, \m{d_0,\ldots,d_{n-1}} des entiers, avec
\m{r>\sigg_{j=0}^{n-1}r_j} et \m{r_i\geq 0} pour \m{0\leq i<n}.

1 - Soit $\kp$ l'ensemble des classes d'isomorphisme de
faisceaux quasi-localement libres sur \m{C_n} de type complet
\m{\big((r_0,\ldots,r_{n-1}),(d_0,\ldots,d_{n-1})\big)} (cf.
\ref{type_complet}). Alors $\kp$ est irr\'eductible.

2 - Soit $\ks$ l'ensemble des classes d'isomorphisme de
morphismes surjectifs \ \m{\E\to F}, o\`u $\E$ est un fibr\'e vectoriel
alg\'ebrique de rang $r$ et de degr\'e $d$ sur \m{C_n} et \ \m{F\in\kp}. Alors
$\ks$ est irr\'eductible.
\end{subsub}

\begin{proof} Elle comporte deux \'etapes. On montre d'abord que
l'irr\'eductibilit\'e de $\kp$ entraine celle de $\ks$. Puis on d\'emontre le
th\'eor\`eme par r\'ecurrence sur $n$.

\bigskip

\textsc{\'Etape 1}

Supposons que $\kp$ soit irr\'eductible. 

Soient \m{f_0:\E_0\to
F_0}, \m{f_1:\E_1\to F_1} des morphismes de $\ks$. Puisque $\kp$ est
irr\'eductible il existe une famille plate $\kf$ de faisceaux quasi localement
libres sur \m{C_n} param\'etr\'ee par une vari\'et\'e alg\'ebrique lisse
irr\'eductible $S$, telle que pour tout point ferm\'e $s$ de $S$, \m{\kf_s} soit
dans \m{\kp} et qu'il existe des points ferm\'es \m{s_0,s_1} de $S$ tels que
\m{\kf_{s_0}\simeq F_0}, \m{\kf_{s_1}\simeq F_1}.

Soit \m{\ko(1)} un fibr\'e en droites tr\`es ample sur \m{C_n}. Si \m{m\gg 0}
les fibr\'es \m{\E_i(m)} sont engendr\'es par leurs sections. Il existe donc
des morphismes surjectifs
\[\phi_i:\ko(-m)\ot\C^{r+1}\lra\E_i .\]
 Posons \m{\D_i=\ker(\phi_i)}, c'est un fibr\'e en droites sur \m{C_n}. On a \
\m{\D_i\simeq\det(\E_i)^{-1}\ot\ko(-(r+1)m)}. Soit \m{N_i} le noyau du
morphisme compos\'e
\xmat{\ko(-m)\ot\C^{r+1}\ar[r]^-{\phi_i} & \E_i\ar[r]^-{f_i} & F_i \ .}

Montrons que \ \m{h^1(N_i\ot\D_i^{-1})=0} \ si \m{m\gg 0}. 
D'apr\`es la proposition \ref{perio_pro1} il existe un fibr\'e vectoriel
\m{\F_i} sur \m{C_n} et un morphisme surjectif
\[\alpha_i:\F_i\lra F_i\]
induisant un isomorphisme \m{\F_{i\mid C}\simeq F_{i\mid C}}. Soit \
\m{A_i=\ker(\alpha_i)}. Si \m{m\gg 0} on a \hfil\break
\m{\Ext^1_{\ko_n}(\ko(-m),A_i)=\nsp}, donc il existe un morphisme
\begin{equation}\label{equ000}\beta_i:\ko(-m)\ot\C^{r+1}\lra\F_i\end{equation}
tel que \ \m{\alpha_i\circ\beta_i=f_i\circ\phi_i}, et \m{\beta_i} est surjectif
(d'apr\`es la proposition \ref{inj_surj}). Soit
\[s=\rg(\F_i)=\rg(F_i)=\sigg_{j=0}^{n-1}r_j .\]
Si \m{m\gg 0}, \m{\F_i(m)} est engendr\'e par ses sections, et \m{s+1}
parmi celles de (\ref{equ000}) suffisent \`a l'engendrer. Il existe donc une
suite exacte
\[0\lra\B_i\lra\ko(-m)\ot\C^{s+1}\lra\F_i\lra 0 ,\]
o\`u \m{\B_i} est un fibr\'e en droites sur \m{C_n}, \m{\B_i\simeq
\det(\F_i)^{-1}\ot\ko(-(s+1)m)}. On a un diagramme commutatif avec lignes et
colonnes exactes
\xmat{ & 0\ar[d] & 0\ar[d]\\
0\ar[r] & \B_i\ar[r]\ar[d] & \ko(-m)\ot\C^{s+1}\ar[r]\ar[d] & \F_i\ar[r]\fleq[d]
& 0\\
0\ar[r] & \ker(\beta_i)\ar[r]\ar[d] & \ko(-m)\ot\C^{r+1}\ar[r]\ar[d] &
\F_i\ar[r] & 0\\
& \ko(-m)\ot\C^{r-s}\fleq[r]\ar[d] & \ko(-m)\ot\C^{r-s}\ar[d]\\
& 0 & 0}
On a
\[\D_i^{-1}(-m)\simeq\det(\E_i)\ot\ko(rm) ,\quad
\B_i\ot\D_i^{-1}\simeq\det(\E_i)\ot\det(\F_i)^{-1}\ot\ko((r-s)m) ,\]
et \m{r-s>0}, donc \ \m{h^1(\D_i^{-1}(-m))=h^1(\B_i\ot\D^{-1})=0} \ si \m{m\gg
0}. On a donc \hfil\break \m{h^1(\ker(\beta_i)\ot\D_i^{-1})=0} \ d'apr\`es la
colonne exacte de gauche du diagramme pr\'ec\'edent. On utilise maintenant la
suite exacte
\[0\lra\ker(\beta_i)\lra N_i\lra\ker(\alpha_i)\lra 0 .\]
Le faisceau \m{\ker(\alpha_i)} ne d\'epend pas de $m$, donc \
\m{h^1(\ker(\alpha_i)\ot\D_i^{-1})=0} \ si \m{m\gg 0}, d'o\`u \
\m{h^1(N_i\ot\D_i^{-1})=0} .

On construit maintenant une famille de morphismes de $\ks$ contenant \m{f_0} et
\m{f_1}. On peut supposer que $m$ est assez grand pour que \m{h^1(\kf_s(m))=0}
pour tout \m{s\in S}.

Soit $U$ la vari\'et\'e des morphismes surjectifs
\[\ko(-m)\ot\C^{r+1}\lra\kf_s ,\]
$s$ parcourant $S$. C'est un ouvert d'un fibr\'e vectoriel sur $S$ (car
\m{h^0(\kf_s(m))} est ind\'ependant de $s$). Il existe un morphisme
universel surjectif de fibr\'es sur \m{U\times C_n}
\[\rho:p_{C_n}^*(\ko(-m))\ot\C^{r+1}\lra(\pi\times I_{C_n})^*(\kf) ,\]
($\pi$ d\'esignant la projection \m{U\to S}). Soit \m{\kn=\ker(\rho)}, c'est une
famille plate de faisceaux quasi localement libres sur \m{C_n} param\'etr\'ee
par $U$.

Soit \ \m{k=\Deg(\D_i)=-d-(r+1)m\Deg(\ko(-1))}. On utilise ici la vari\'et\'e
\m{\Lambda^k(C_n)} et le ``fibr\'e de Poincar\'e'' $\kd$ sur \
\m{\Lambda^k(C_n)\times C_n} (cf. \ref{pic_cn}). Soit $V$ la vari\'et\'e des
morphismes \m{\kd_y\to\kn_u}, \m{y\in\Lambda^k(C_n)}, \m{u\in U}, tels que le
compos\'e \ \m{\kd_y\to\kn_u\subset\ko(-m)\ot\C^{r+1}} \ soit un morphisme
injectif de fibr\'es vectoriels, et que \ \m{h^1(\kd_y^{-1}\ot\kn_u)=0} . C'est
un ouvert d'un fibr\'e vectoriel sur \ \m{U\times\Lambda^k(C_n)}. On a un
morphisme universel
\[\Theta:p_{\Lambda^k(C_n)\times C_n}^*(\kd)\lra p_{C_n}^*(\ko(-m))\ot\C^{r+1}
.\]
Soit \m{\E=\coker(\Theta)}, c'est un fibr\'e vectoriel sur \m{V\times C_n}. On
a un morphisme canonique surjectif
\[\E\lra (\tau\times I_{C_n})^*(\kf) ,\]
$\tau$ d\'esignant la suite de projections \m{V\to U\to S}. C'est la famille de
morphismes de $\kp$ contenant \m{f_0} et \m{f_1}.

\bigskip

\textsc{\'Etape 2}

On montre maintenant que si le th\'eor\`eme est vrai pour \m{n-1},
c'est-\`a-dire sur \m{C_n}, alors $\kp$ est irr\'eductible.
Soient \m{F_0,F_1\in\kp}. D'apr\`es la proposition \ref{perio_pro1} il existe
des fibr\'es vectoriels \m{\F_i} sur \m{C_n} et des morphismes surjectifs \
\m{\phi_i:\F_i\to F_i} \ induisant des isomorphismes \ \m{\F_{i\mid C}\simeq
F_{i\mid C}}. Les faisceaux \m{\kn_i=\ker(\phi_i)} sont quasi localement libres
de m\^eme type complet, d'apr\`es le corollaire \ref{perio_cor}, et leur
support est contenu dans \m{C_{n-1}}. On consid\`ere maintenant les morphismes
surjectifs de faisceaux sur \m{C_{n-1}}, induits par les transpos\'es des
\m{\phi_i}
\[g_i:\F_{i\mid C_{n-1}}^\vee\lra\kn_i^\vee .\]
D'apr\`es l'hypoth\`ese de r\'ecurrence il existe des familles plates $\G$,
$\ke$ de faisceaux coh\'erents sur \m{C_{n-1}} param\'etr\'ees par une
vari\'et\'e alg\'ebrique lisse irr\'eductible $S$, et un morphisme surjectif
\m{\gamma:\G\to\kf} tels que :
\begin{enumerate}
\item[--] $\G$ est localement libre.
\item[--] Les faisceaux $\kf_s$, $s\in S$, sont de m\^eme type complet que
$\kn_0^\vee$, $\kn_1^\vee$.
\item[--] Il existe des points ferm\'es $s_0,s_1\in S$ tels que \
$\gamma_{s_0}\simeq g_0$ \ et \ $\gamma_{s_1}\simeq g_1$.
\end{enumerate}
D'apr\`es le th\'eor\`eme \ref{prol_th} il existe un prolongement de $\G$ en un
fibr\'e vectoriel $\E$ sur \m{C_n} tel que \m{\E_{s_i}=\F_i^\vee}. On
consid\`ere maintenant le morphisme compos\'e surjectif de faisceaux sur \m{C_n}
\xmat{\Theta:\E\flon[r] & \G\ar[r]^-\gamma & \kn^\vee .}
Soit \ \m{\kv=\ker(\Theta)}. C'est une famille plate de faisceaux
quasi-localement libres sur \m{C_n} de m\^eme type complet que
\m{F_0^\vee,F_1^\vee} (d'apr\`es la remarque \ref{perio_rem}). La famille
contenant \m{F_0,F_1} que l'on cherche est alors \m{\kv^\vee}.
\end{proof}

\end{sub}

\sepsub

\Ssect{Faisceaux quasi localement libres de type rigide}{rig-qll}

Soit $\ke$ un faisceau coh\'erent quasi localement libre sur \m{C_n}. Soient
\m{a=\lbrack\frac{R(\ke)}{n}\rbrack} et \m{k=R(\ke)-an}. On a donc
\m{R(\ke)=an+k}. On dit que $\ke$ est de {\em type rigide} s'il est localement
libre si \m{k=0} et localement isomorphe \`a \m{a\ko_n\oplus\ko_k} si \m{k>0}.
Si \m{k>0} cela revient � dire que $\ke$ est de type \m{(m_1,\ldots,m_n)}, avec
\m{m_i=0} si \m{i\not=k,n} et \m{m_k=0} ou $1$.

\sepprop

\begin{subsub}\label{qlltr1}{\bf Proposition : } Soient $Y$ une vari\'et\'e
alg\'ebrique int\`egre et $\kf$ une famille plate de faisceaux coh\'erents sur
\m{C_n} param\'etr\'ee par $Y$. Alors l'ensemble des points \m{y\in Y} tels que
\m{\ke_y} soit quasi localement libre de type rigide est un ouvert de $Y$.
\end{subsub}

\begin{proof} Soit $R$ le rang g\'en\'eralis\'e des faisceaux \m{\ke_y},
\m{R=an+k}, avec \m{a\geq 0}, \m{0\leq k<n}. On supposera que \m{k>0} car le
cas \m{k=0} (o\`u quasi localement libre de type rigide \'equivaut \`a
localement libre) est bien connu.

Soit \m{y_0\in Y} tel que \m{\ke_{y_0}} soit quasi localement
libre de type rigide. On va montrer que \m{\ke_y} l'est aussi pour tout $y$
dans un voisinage de \m{y_0}.

Soient \m{y\in Y} et \m{P\in C} tels qu'il existe un morphisme surjectif \
\m{\phi:(a+1)\ko_{nP}\to\ke_{yP}}. On peut prolonger $\phi$ en un
morphisme surjectif \m{(a+1)\ko_U\to\ke_{\mid U}}, $U$ d\'esignant un voisinage
convenable de \m{(y,P)}. Les points \m{(y,P)} poss\'edant cette propri\'et\'e
constituent donc un ouvert $W$ de \m{Y\times C_n}.

Posons \m{M=\ke_{yP}}. Soit 
\[M_n=\nsp\subset M_{n-1}\subset\cdots\subset M_{1}\subset M_0=M \]
la premi\`ere filtration canonique de $M$. Alors on a \ \m{\rg(M/M_1)\geq a+1} ,
car dans le cas contraire on aurait \ \m{R(M)\leq an}. Puisque $\phi$ induit un
morphisme surjectif \ \m{\psi:(a+1)\ko_{CP}\to M/M_1} , on a \
\m{\rg(M/M_1)=a+1}, et $\psi$ est un isomorphisme. On a alors \
\m{\rg(M_i/M_{i+1})=a+1} \ pour \m{0\leq i<k} et \ \m{\rg(M_i/M_{i+1})=a} \
pour \m{k\leq i<n}. Puisque la multiplication par une \'equation de $C$ induit
des morphismes surjectifs \ \m{M/M_1\to M_i/M_{i+1}} \ on a \ \m{M_i/M_{i+1}
\simeq (a+1)\ko_{nP}}. Donc \m{\ke_{yP}/(\ke_{yP})_k} est un
\m{\ko_{kP}}-module libre.

D'apr\`es le corollaire \ref{f_ql_cor1}, en on point \m{(y,P)} de
$W$, la fibre de \m{(\ke_{\mid W})_k} n'est autre que \m{(\ke_{yP})_k}. Soit
\m{V\subset W} l'ouvert constitu\'e des points \m{(y,P)} tels que
\m{(\ke_{yP})_k} soit un \m{\ko_{P,n-k}}-module libre. Soit \m{Z=(Y\times
C_n)\backslash V}, et $T$ sa projection sur $Y$, qui est une sous-vari\'et\'e
ferm\'ee. L'ouvert \m{Y\backslash T} est le voisinage recherch\'e de \m{y_0}.
\end{proof}

\sepprop

\begin{subsub}D\'eformations des faisceaux quasi localement libres de type
rigide - \rm Soit $\ke$ un faisceau coh\'erent sur \m{C_n}. Soit
\m{(\widetilde{\ke},S,s_0,\epsilon)} une {\em d\'eformation semi-universelle}
de $\ke$ (cf. \cite{si_tr}, \cite{dr3} 3.1), donc $\widetilde{\ke}$ est une
famille plate de faisceaux coh\'erents sur \m{C_n} param\'etr\'ee par $S$,
\m{s_0} est un point ferm\'e de $S$ et \
\m{\epsilon:\widetilde{\ke}_{s_0}\simeq\ke}. Le morphisme de d\'eformation
infinit\'esimale de Koda\"ira-Spencer
\[\omega_{\widetilde{\ke},s_0}:T_{s_0}S\lra\Ext^1_{\ko_n}(\ke,\ke)\]
est un isomorphisme. On pose
\[D_{reg}(\ke) \ = \ \omega_{\widetilde{\ke},s_0}(T_{s_0}(S_{red})) .\]
Si $\kf$ est une d\'eformation de $\ke$ param\'etr\'ee par une vari\'et\'e
alg\'ebrique r\'eduite $Y$, et si \m{\kf_y\simeq\ke}, l'image du morphisme de
d\'eformation infinit\'esimale de Koda\"ira-Spencer \m{\omega_{\kf,y}} est
contenue dans \m{D_{reg}(\ke)}.

On dit que $\ke$ est {\em lisse pour les d\'eformations r\'eduites} si
\m{S_{red}} est lisse en \m{s_0}.

Si $\ke$ est un faisceau quasi localement libre de type rigide sur \m{C_n}, on a
\[D_{reg}(\ke) \ \subset \ H^1(\EEnd(\ke))\]
d'apr\`es la proposition \ref{qlltr1} et le corollaire \ref{f_ql_cor3}.
\end{subsub}

\sepprop

\begin{subsub}\label{qlltr2}{\bf Th\'eor\`eme : } Si $\ke$ est un faisceau quasi
localement libre de type rigide g\'en\'erique, alors on a \ \m{D_{reg}(\ke)=
H^1(\EEnd(\ke))}.
\end{subsub}

\begin{proof} Par r\'ecurrence sur $n$. Le r\'esultat est vrai pour \m{n=1}.
Supposons le vrai pour \m{n-1\geq 1}. On peut donc supposer que \
\m{D_{reg}(\ke_1)= H^1(\EEnd(\ke_1))} \ sur \m{C_{n-1}}. A partir d'une
d\'eformation compl\`ete de $\ke_1$ (comme faisceau sur \m{C_{n-1}}) on va
construire une d\'eformation compl\`ete de $\ke$ en utilisant les r\'esultats
de \ref{constr}, puis montrer que le morphisme de d\'eformation infinit\'esimale
de Koda\"ira-Spencer de cette d\'eformation en $\ke$ a pour image
\m{H^1(\EEnd(\ke))}.

On a une suite exacte
\[0\lra\HHom(\ke_{\mid C},\ke^{(1)})\lra\EEnd(\ke)\lra\EEnd(\ke_1)\lra 0 ,\]
d'o\`u on d\'eduit la suivante
\begin{equation}\label{equ_qlltr0}\Ext^1_{\ko_C}(\ke_{\mid C},\ke^{(1)})\lra
H^1(\EEnd(\ke))\lra H^1(\EEnd(\ke_1))\lra 0 .
\end{equation}
Soit $\kf$ une famille plate de faisceaux quasi localement libres sur
\m{C_{n-1}}, de m\^eme type complet que \m{\ke_1}, param\'etr\'ee par une
vari\'et\'e int\`egre $Y$, telle qu'il existe \m{y_0\in Y} tel que
\m{\kf_{y_0}\simeq\ke_1} et que l'image du morphisme de d\'eformation
infinit\'esimale de Koda\"ira-Spencer \m{\omega_{\kf,y_0}} soit \
\m{H^1(\EEnd(\ke_1))}. On peut supposer que \m{h^1(\EEnd(\ke_1))} est minimal,
donc \m{h^1(\EEnd(\kf_y))} est constant au voisinage de \m{y_0}. En rempla\c
cant $Y$ par un voisinage de \m{y_0} on peut donc supposer que
\begin{enumerate}
\item[--] $h^1(\EEnd(\kf_y)$ est ind\'ependant de $y\in Y$,
\item[--] pour tout $y\in Y$, l'image de \m{\omega_{\kf,y}} est
$H^1(\EEnd(\kf_y)$.
\end{enumerate}
Puisque $\ke$ est quasi localement libre de type rigide on a \
\m{\Gamma^{(0)}(\ke)=0} \ ou bien c'est un fibr\'e en droites, ce dernier cas
ne pouvant se produire que si \m{\ke_1} est localement libre. Dans tous les cas
on a des suites exactes
\begin{equation}\label{equ_qlltr}
0\lra\Gamma^{(0)}(\ke)\lra\ke_{\mid C}\lra\ke_{1\mid C}\ot L^*\lra 0 ,
\end{equation}
\begin{equation}\label{equ_qlltrb}
0\lra(\ke_1)^{(1)}\lra\ke^{(1)}\lra\Gamma^{(0)}(\ke)\lra 0 ,
\end{equation}

d'apr\`es le corollaire \ref{filt3} et la proposition \ref{filt1}, (iv).

\bigskip

\textsc{\'Etape 1 - }
{\em Param\'etrisation des restrictions \`a $C$ des d\'eformations de $\ke$}

On suppose que \m{\Gamma^{(0)}(\ke)\not=0}, soit $d$ son degr\'e. On peut
supposer que \hfil\break \m{\dim_\C(\Ext^1_{\ko_C}(\ke_{1\mid C}\ot L^*,
\Gamma^{(0)}(\ke)))} \ est minimal parmi les \ \m{\dim_\C(\Ext^1_{\ko_C}(
\kf_{y\mid C}\ot L^*,D))} , \m{y\in Y}, \m{D\in\Pic^d(C)}. Soit $W$ l'ouvert
des points \m{(y,D)} de \m{Y\times\Pic^d(C)} tels que \hfil\break
\m{\dim_\C(\Ext^1_{\ko_C}(\kf_{y\mid C}\ot L^*,D))} soit minimal. Soit $\W$ le
fibr\'e relatif
\[\W \ = \EExt^1_p(p_Y^\sharp(\kf_{\mid C}\ot L^*),p_d^\sharp(\kd)) ,\]
$p$, \m{p_Y}, \m{p_d} d\'esignant les projections \ \m{Y\times\Pic^d(C)\times
C\to Y\times\Pic^d(C)} , \m{Y\times\Pic^d(C)\to Y}, \m{Y\times\Pic^d(C)\to
\Pic^d(C)} respectivement, et $\kd$ un fibr\'e de Poincar\'e sur \m{\Pic^d(C)
\times C}. Il existe une {\em extension universelle} sur \m{\W\times C}
\xmat{0\ar[r] & {p'_d}^\sharp(\kd)\ar[r] &
\kv\ar[r]^-\theta & {p'_Y}^\sharp(\kf_{\mid C}\ot L^*)\ar[r] &  0 ,}
\m{p'_d}, \m{p'_Y} d\'esignant les projections \m{\W\to Pic^d(C)}, \m{\W\to Y}
respectivement.

Si \m{\Gamma^{(0)}(\ke)=0} on prend \m{\W=Y} et \m{\kv=\kf_{\mid C}\ot L^*},
\m{p'_Y} d\'esigne l'identit\'e \m{\W\to Y} et $\theta$ l'identit\'e \m{\kv\to
p_Y^\sharp(\kf_{\mid C}\ot L^*)}.

On note \m{w_0} le point de $\W$ correspondant \`a $\ke$ (plus exactement
correspondant \`a l'extension (\ref{equ_qlltr})). On a \m{p'_Y(w_0)=y_0}.

\bigskip

\textsc{\'Etape 2 - }
{\em Param\'etrisation des d\'eformations de $\ke$}

On peut supposer que \m{\dim_\C(\Ext^1_{\ko_C}(\ke_{\mid C},(\ke_1)^{(1)}))} est
minimal parmi les \hfil\break\m{\dim_\C(\Ext^1_{\ko_C}(\kv_w,(\kf_y)^{(1)}))},
\m{w\in\W}, \m{y=p'_Y(w)}. On peut supposer aussi que \hfil\break
\m{\dim_\C(\Ext^1_{\ko_n}( \ke_C,\ke_1))} est minimal parmi les
\m{\dim_\C(\Ext^1_{\ko_n}(\kv_w,\kf_y))}. Soit $U$ l'ouvert de $\W$ o\`u ces
deux dimensions sont minimales. On consid\`ere le fibr\'e relatif
\[\U \ = \ \EExt^1_{p_U}(\kv,{p'_Y}^\sharp(\kf)) ,\]
\m{p_U} d\'esignant la projection \m{U\times C_n\to U}. On a un morphisme
canonique
\[\tau : \U\lra p_{U*}(\HHom(\kv\ot L,{p'_Y}^\sharp(\kf_{\mid C}))\]
qui en \m{w\in U}, \m{y={p_Y'}(w)}, est le morphisme canonique
\xmat{\Ext^1_{\ko_n}(\kv_w,\kf_y)\ar[r] & \Hom(\kv_w\ot L,\kf_{y\mid C})}
(cf. \ref{constr}). On a une section canonique de \m{{p'_Y}^\sharp(\kf_{\mid
C})} induite par $\theta$. Soit \m{\T\subset{p'_Y}^\sharp(\kf_{\mid
C})} la sous-vari\'et\'e correspondante (isomorphe \`a $U$).

La vari\'et\'e \ \m{Z=\tau^{-1}(\T)} \ est un fibr\'e en espaces affines. On a
une extension universelle sur \m{Z\times C_n}
\[0\lra q_Y^\sharp(\kf)\lra\E\lra q_U^\sharp(\kv)\lra 0 ,\]
\m{q_Y}, \m{q_U} d\'esignant les projections \m{Z\to Y} et \m{Z\to U}
respectivement. Soient \m{\sigma\in Z} au dessus de \m{w\in U}, et
\m{y=p'_Y(w)}. Alors, sur \m{\lbrace\sigma\rbrace\times C_n}, l'extension
pr\'ec\'edente est
\[0\lra\kf_y\lra\E_\sigma\lra\kv_w\lra 0\]
associ\'ee \`a \m{\sigma\in\Ext^1_{\ko_n}(\kv_w,\kf_y)}. L'image de $\sigma$
dans \m{\Hom(\kv_w\ot L,\kf_{y\mid C})} est \m{\theta_w}.

On note \m{z_0} le point de $Z$ correspondant \`a $\ke$ (ou plus exactement \`a
l'extension \hfil\break \m{0\to\ke_1\to\ke\to\ke_{\mid C}\to 0}) .

\bigskip

\textsc{\'Etape 3 - }{\em Surjectivit\'e de \m{\omega_{\E,z_0}}}

On a un diagramme commutatif
\xmat{T_{e_0}Z\flon[r]^-{Tq_Y}\ar[d]^{\omega_{\E,e_0}} &
T_{y_0}Y\flon[d]^{\omega_{\kf,y_0}} \\ H^1(\EEnd(\ke))\flon[r]^-\rho &
H^1(\EEnd(\ke_1))}
Il en d\'ecoule que \m{\rho\omega_{\E,e_0}} est surjectif. Il reste \`a voir
que le morphisme induit par \m{\omega_{\E,e_0}}
\[\nu : \ker(Tq_Y)\lra\ker(\rho)\]
est surjectif.

Supposons d'abord que \m{\Gamma^{(0)}(\ke)=0}. Dans ce cas on a
\[\ker(Tq_Y) \ = \ \Ext^1_{\ko_C}(\ke_{\mid C},\ke^{(1)}) \]
et un morphisme surjectif
\[\Ext^1_{\ko_C}(\ke_{\mid C},\ke^{(1)})\lra\ker(\rho)\]
d'apr\`es (\ref{equ_qlltr0}). Il suffit donc de montrer que le diagramme
\xmat{\Ext^1_{\ko_C}(\ke_{\mid C},\ke^{(1)})\flinc[r]\fleq[d] &
T_{e_0}Z\ar[d]^{\omega_{\E,e_0}} \\
\Ext^1_{\ko_C}(\ke_{\mid C},\ke^{(1)})\flon[r] & H^1(\EEnd(\ke))}
est commutatif. Ce sera fait dans l'\'etape 4.

Supposons maintenant que \m{\Gamma^{(0)}(\ke)\not=0}. Soient \
\m{A={p_Y'}^{-1}(y_0)\subset\W} \ et \m{w_0\in A} correspondant \`a
\m{\ke_{\mid C}} (on a donc \m{p'_Y(w_0)=y_0}). 

Soit
\xmat{\beta:\Ext^1_{\ko_C}(\ke_{\mid C},(\ke_1)^{(1)})\ar[r] &
\Ext^1_{\ko_C}(\ke_{\mid C},\ke_1)\flon[r] & \ker(\rho)}
le morphisme compos\'e (cf. (\ref{equ_qlltr0})). D'apr\`es la suite
exacte (\ref{equ_qlltrb}) on a un morphisme surjectif
\[\psi:\Ext^1_{\ko_C}(\ke_{\mid
C},\Gamma^{(0)}(\ke))\lra\ker(\rho)/\imm(\beta) .\]

Si \m{(\widetilde{\ke},S,s_0,\epsilon)} est une d\'eformation semi-universelle
de $\ke$, on a, en consid\'erant le fibr\'e en droites
\m{\Gamma^{(0)}(\widetilde{\ke})} sur \m{S\times C} (cf. corollaire
\ref{f_ql_cor1}) un morphisme de d\'eformation infinit\'esimale de
Koda\"ira-Spencer
\[\Ext^1_{\ko_n}(\ke,\ke)\lra\Ext^1_{\ko_C}(\Gamma^{(0)}(\ke),\Gamma^{(0)}
(\ke))\]
qui s'annule sur \m{\imm(\beta)}. On obtient donc un morphisme canonique
\[\omega_0:\ker(\rho)/\imm(\beta)\lra\Ext^1_{\ko_C}(\Gamma^{(0)}(\ke),\Gamma^{
(0)}(\ke)) .\]
On a des suites exactes
\xmat{0\ar[r] & \Ext^1_{\ko_C}(\ke_{\mid C},(\ke_1)^{(1)})\ar[r]^-\gamma &
\ker(Tq_Y)=T_\sigma(q_Y^{-1}(y_0))\ar[r] & T_{w_0}A\ar[r] & 0 ,}
\[0\lra\Ext^1_{\ko_C}(\ke_{1\mid C}\ot L^*,\Gamma^{(0)}(\ke))\lra T_{w_0}A\lra
\Ext^1_{\ko_C}(\Gamma^{(0)}(\ke),\Gamma^{(0)}(\ke))\lra 0 ,\]
et comme dans le cas \m{\Gamma^{(0)}(\ke)=0} un diagramme commutatif
\begin{equation}\label{equ_qlltr1}
\xymatrix{\Ext^1_{\ko_C}(\ke_{\mid C},(\ke_1)^{(1)})\fleq[d]\ar[r]^-\gamma &
\ker(Tq_Y)\ar[d]^{\omega_{\E,e_0}} \\ \Ext^1_{\ko_C}(\ke_{\mid C},(\ke_1)^{(1)})
\ar[r]^-\beta & \ker(\rho).}
\end{equation}
Il suffit donc de montrer la surjectivit\'e du morphisme \
\m{\theta:T_{w_0}A\to\ker(\rho)/\imm(\beta)} \ induit par \m{\omega_{\E,e_0}}.

On a un diagramme commutatif
\xmat{\Ext^1_{\ko_C}(\ke_{1\mid C}\ot L^*,\Gamma^{(0)}(\ke))\fleq[dd]\flinc[r] &
T_{w_0}A\ar[d]^\theta\flon[r] & \Ext^1_{\ko_C}(\Gamma^{(0)}(\ke),\Gamma^{
(0)}(\ke))\fleq[dd] \\
& \ker(\rho)/\imm(\beta)\ar[ur]^{\omega_0}\\
\Ext^1_{\ko_C}(\ke_{1\mid C}\ot L^*,\Gamma^{(0)}(\ke))\ar[r]
& \Ext^1_{\ko_C}(\ke_{\mid C},\Gamma^{(0)}(\ke))\flon[u]^\psi\flon[r] &
\Ext^1_{\ko_C}(\Gamma^{(0)}(\ke),\Gamma^{(0)}(\ke))
}
(voir \'etape 4). La surjectivit\'e de $\theta$ en d\'ecoule ais\'ement.

\bigskip

\textsc{\'Etape 4 - }{\em Commutativit\'e des diagrammes}

On ne d\'emontrera que la commutativit\'e du diagramme (\ref{equ_qlltr1}), les
autres cas \'etant analogues. Rappelons que
\m{z_0\in\Ext^1(\ke_{\mid C},\ke_1)} d\'esigne l'\'el\'ement associ\'e \`a
l'extension \m{0\to\ke_1\to\ke\to\ke_{\mid C}\to 0}. Soit
\[\alpha\in\Ext^1_{\ko_C}(\ke_{\mid C},(\ke_1)^{(1)})\subset\Ext^1(\ke_{\mid
C},\ke_1) .\]
Soit $\kg$ la famille plate de faisceaux quasi localement libres sur \m{C_n}
param\'etr\'ee par $\C$ d\'efinie par~: pour tout \m{\lambda\in\C},
\m{0\to\ke_1\to\kg_\lambda\to\ke_{\mid C}\to 0} est l'extension associ\'ee \`a
\m{z_0+\lambda\alpha}. On consid\`ere le morphisme de d\'eformation
infinit\'esimale de Koada\"ira-Spencer de cette famille en \m{\lambda=0} :
\[\omega_{\kg,0}:\C\lra H^0(\EEnd(\ke)) .\]
Il suffit de montrer que \m{\omega_{\kg,0}(1)} est \'egal \`a l'image de
$\alpha$ dans \m{H^0(\EEnd(\ke))}. On utilisera les r\'esultats de \ref{extens}.
Soient \m{\phi:\spec(\C[t]/(t^2))\to\C} le morphisme induit par le quotient
\m{\C[t]\to\C[t]/(t^2)} et \ \m{\pi:\spec(\C[t]/(t^2))\times C_n\to C_n} \ la
projection. Alors le faisceau \m{\pi_*(\phi^\sharp(\kg))} est une extension de
$\ke$ par lui-m\^eme, et \m{\omega_{\kg,0}(1)} n'est autre que l'\'el\'ement de
\m{\Ext^1_{\ko_n}(\ke,\ke)} associ\'e \`a cette extension.(cf. \cite{dr3}, 3-).

Il existe un recouvrement ouvert \m{(U_i)} de $C$ tel que $\alpha$ soit
repr\'esent\'e par un cocycle \m{(\alpha_{ij})}, \m{\alpha_{ij}:
\ke_{\mid C\mid U_{ij}}\to(\ke_1)^{(1)}_{\mid U_{ij}}}. D'apr\`es la
proposition \ref{ext_P}, \m{\kg_\lambda}
est obtenu en recollant les \m{\ke_{\mid U_{ij}}} au moyen des automorphismes
\m{I+\lambda\tau_{ij}} de \m{\ke_{\mid U_{ij}}}, \m{\tau_{ij}} \'etant le
compos\'e
\xmat{\ke_{\mid U_{ij}}\flon[r] & \ke_{\mid C\mid U_{ij}}\ar[r]^-{\alpha_{ij}} &
(\ke_1)^{(1)}_{\mid U_{ij}}\flinc[r] & \ke_{1\mid U_{ij}}\flinc[r] &
\ke_{\mid U_{ij}} .}
Donc $\kg$ est obtenu en recollant les \m{\pi^*(\ke)_{\mid\C\times U_{ij}}} au
moyen des automorphismes \m{I+t\tau_{ij}} (avec \m{t=I_\C}). Il en d\'ecoule
que \m{\pi_*(\phi^\sharp(\kg))} est obtenu en recollant les
\m{(\ke\oplus\ke)_{\mid U_{ij}}} au moyen des automorphismes \m{\begin{pmatrix}
1 & \alpha_{ij}\\ 0 & 1\end{pmatrix}}. D'apr\`es la proposition \ref{ext_P},
l'\'el\'ement de \m{\Ext^1_{\ko_n}(\ke,\ke)} associ\'e \`a l'extension
de $\ke$ par lui-m\^eme obtenue ainsi n'est autre que l'image de $\alpha$.
\end{proof}

\sepprop

\begin{subsub}\label{qlltr3}{\bf Corollaire : } Si $\ke$ est un faisceau
quasi localement libre de type rigide tel que \m{\dim_\C(\End(\ke))} soit
minimal, c'est \`a-dire soit tel que pour tout faisceau quasi localement libre
$\kf$ de m\^eme type complet que $\ke$ on ait \ \m{\dim_\C(\End(\kf))\geq
\dim_\C(\End(\ke))}, alors on a \ \m{D_{reg}(\ke)= H^1(\EEnd(\ke))} , et $\ke$
est lisse pour les d\'eformations r\'eduites.
\end{subsub}

\begin{proof} Soit \m{d=\dim_\C(H^1(\EEnd(\ke)))}. Soit
\m{(\widetilde{\ke},S,s_0,\epsilon)} une d\'eformation semi-universelle de
$\ke$. Alors, si $ \ke$ est g\'en\'erique, d'apr\`es le th\'eor\`eme
\ref{qlltr2}, l'espace tangent de \m{S_{red}} en \m{s_0} est \m{H^1(\EEnd(\ke))}
et $S$ est lisse en \m{s_0}. Donc dans tous les cas la dimension de \m{S_{red}}
est aussi \'egale \`a $d$. Puisque l'espace tangent \m{T_{s_0}S_{red}} est de
dimension au moins $d$ et est contenu dans \m{H^1(\EEnd(\ke))}, il lui est
\'egal et le corollaire \ref{qlltr3} en d\'ecoule.
\end{proof}
\end{sub}

\sepsub

\Ssect{Faisceaux quasi localement libres de type rigide stables}{f_st}

La notion de (semi-)stabilit\'e des faisceaux coh\'erents sur \m{C_n} est
ind\'ependante du choix d'un fibr\'e en droites ample sur \m{C_n} (cf.
\cite{dr2}). Un faisceau coh\'erent $\ke$ sur \m{C_n} est dit {\em semi-stable}
(resp. stable) s'il est sans torsion et si pour tout sous-faisceau propre
\m{\kf\subset\ke} on a \ \m{\mu(\kf)<\mu(\ke)} \ (resp.
\m{\mu(\kf)\leq\mu(\ke)}) .

Soient $R$, $D$ des entiers, avec \m{R\geq 1}. On note \m{\km(R,D)} la
vari\'et\'e de modules des faisceaux stables de rang g\'en\'eralis\'e $R$
et de degr\'e g\'en\'eralis\'e $D$ sur \m{C_n} (cf. \cite{ma1}, \cite{ma2},
\cite{si}). On supposera que \m{\deg(L)<0}, car dans le cas contraire les
seuls faisceaux sans torsion stables sur \m{C_n} sont les fibr\'es vectoriels
stables sur $C$.

Soit $\ke$ un faisceau quasi localement libre de type rigide, non localement 
libre. Il existe donc des entiers \m{a=a(\ke),k=k(\ke)}  tels que
\m{a,k> 0}, \m{k<n}, et que $\ke$ soit localement isomorphe \`a
\ \m{a\ko_n\oplus\ko_k} . Posons
\[E=\ke_{\mid C} , \quad\quad F=G_k(\ke)\ot L^{-k} .\]
Alors on a \ \m{\rg(E)=a+1} , \m{\rg(F)=a}, et
\[(G_0(\ke),G_1(\ke),\ldots,G_{n-1}(\ke)) \ = \ (E,E\ot L,\ldots,E\ot
L^{k-1},F\ot L^k,\ldots, F\ot L^{n-1}) \ .\]
Donc
\[\Deg(\ke) \ = \ k\deg(E)+(n-k)\deg(F)+\big(n(n-1)a+
k(k-1)\big)\deg(L)/2 \ .\]
Posons \ \m{\delta=\delta(\ke)=\deg(F)}, \m{\epsilon=\epsilon(\ke)=\deg(E)} .
D'apr\`es la proposition \ref{qlltr1} les d\'eformations de $\ke$ sont des
faisceaux quasi localement libres de type rigide, et d'apr\`es le corollaire
\ref{f_ql_cor1}, \m{a(\ke)}, \m{k(\ke)}, \m{\delta(\ke)} et \m{\epsilon(\ke)}
sont aussi invariants par d\'eformation. Soient
\[R=an+k \ , \quad\quad D=k\epsilon+(n-k)\delta+
\big(n(n-1)a+k(k-1)\big)\deg(L)/2 \ .\]
Les faisceaux quasi localement libres de type rigide stables $\kf$ tels que
\m{a(\kf)=a}, \m{k(\kf)=k}, \m{\delta(\kf)=\delta}, \m{\epsilon(\kf)=\epsilon}
constituent donc un ouvert de \m{\km(R,D)}, not\'e \m{\kn(a,k,\delta,\epsilon)}.

\sepprop

\begin{subsub}\label{f_st_1}{\bf Proposition : } La vari\'et\'e
\m{\kn(a,k,\delta,\epsilon)} est irr\'eductible, et la sous-vari\'et\'e
r\'eduite sous-jacente \m{\kn(a,k,\delta,\epsilon)_{red}} est lisse. Si cette
vari\'et\'e est non vide, on a
\[\dim(\kn(a,k,\delta,\epsilon)) \ = \ 1 - \big(\frac{n(n-1)}{2}a^2+k(n-1)a+
\frac{k(k-1)}{2}\big)\deg(L)+(g-1)(na^2+k(2a+1))\]
($g$ d\'esignant le genre de $C$). Pour tout faisceau $\kf$ de
\m{\kn(a,k,\delta,\epsilon)_{red}}, l'espace tangent de
\m{\kn(a,k,\delta,\epsilon)_{red}} en $\kf$ est canoniquement isomorphe \`a
\m{H^1(\EEnd(\kf))} .
\end{subsub}

\begin{proof} Cela d\'ecoule du fait que tout faisceau stable sur \m{C_n} est
simple, du corollaire \ref{qlltr3} et de la proposition \ref{f_end}.
\end{proof}

\end{sub}

\end{document}